\documentclass[preprint,12pt]{elsarticle}

\usepackage{amssymb}
\usepackage{amsmath}
\usepackage{bm}
\usepackage{color}
\usepackage{exscale}
\usepackage{relsize}
\usepackage{graphicx}
\usepackage{wrapfig}
\usepackage{multicol}
\usepackage{verbatim}
\newtheorem{thm}{Theorem}

\newtheorem{lem}{Lemma}

\begin{document}

\begin{frontmatter}



\title{Higher-order multi-scale deep Ritz method for multi-scale problems of authentic composite materials}

\author[label1]{Jiale Linghu}
\author[label1]{Hao Dong\corref{cor1}}\ead{donghao@mail.nwpu.edu.cn}
\cortext[cor1]{Corresponding author.}
\author[label2]{Junzhi Cui}
\author[label3]{Yufeng Nie}
\address[label1]{School of Mathematics and Statistics, Xidian University, Xi'an 710071, PR China}
\address[label2]{Academy of Mathematics and Systems Science, Chinese Academy of Sciences, Beijing 100190, PR China}
\address[label3]{School of Mathematics and Statistics, Northwestern Polytechnical University, Xi'an 710129, PR China}

\begin{abstract}
The direct deep learning simulation for multi-scale problems remains a challenging issue. In this work, a novel higher-order multi-scale deep Ritz method (HOMS-DRM) is developed for thermal transfer equation of authentic composite materials with highly oscillatory and discontinuous coefficients. In this novel HOMS-DRM, higher-order multi-scale analysis and modeling are first employed to overcome limitations of prohibitive computation and Frequency Principle when direct deep learning simulation. Then, improved deep Ritz method are designed to high-accuracy and mesh-free simulation for macroscopic homogenized equation without multi-scale property and microscopic lower-order and higher-order cell problems with highly discontinuous coefficients. Moreover, the theoretical convergence of the proposed HOMS-DRM is rigorously demonstrated under appropriate assumptions. Finally, extensive numerical experiments are presented to show the computational accuracy of the proposed HOMS-DRM. This study offers a robust and high-accuracy multi-scale deep learning framework that enables the effective simulation and analysis of multi-scale problems of authentic composite materials.
\end{abstract}


\end{frontmatter}


\section{Introduction}
With the rapid development of material science, composite materials are widely used in a variety of engineering sectors, such as aviation, aerospace and civil engineering due to their excellent physical performance. Composite materials possess highly heterogeneous properties, which arise significantly multi-scale characteristics in their spatial structures. It is well-known that classical numerical methods (including finite difference method, finite element method, finite volume method and boundary element method, etc) have been employed to simulate multi-scale problems of composite materials. However, it is extremely difficult and challenge to attain full-scale and high-resolution simulation for large-scale composites via above-mentioned traditional numerical approaches since the decisive disadvantage is the tremendous computational effort, which results from the numerical solution of multi-scale problems of composite materials with fine computation mesh \cite{R1,R45,R46}. Therefore, there is an urgent need to develop novel computational methods that can efficiently simulate multi-scale problems of the composites.

To circumvent time-consuming task and economize computational resources, in the past decades, a considerable number of multi-scale computational methods in research and engineering fields have devised and focused attention on computing and analyzing the physical properties of composite materials, such as asymptotic homogenization method (AHM) \cite{R2,R7,R12}, numerical homogenization method (NHM) \cite{R3}, heterogeneous multi-scale method (HMM) \cite{R4}, variational multi-scale method (VMS) \cite{R5} and multi-scale finite element method (MsFEM) \cite{R6,R13}, etc. In recent years, Cui and his research team systematically developed a class of higher-order multi-scale methods \cite{R8,R9,R10,R11}. By establishing novel higher-order corrected terms in classical asymptotic homogenization method, the numerical accuracy of higher-order multi-scale methods is significantly improved for simulating authentic composite materials, which have been expanded to  multi-physics coupling problems, stochastic multi-scale problem, nonlinear multi-scale problem and structural mechanics problem of heterogeneous materials \cite{R8,R38,R39,R40}. Nevertheless, it should be underlined that above-mentioned multi-scale numerical approaches are still not easy or even ineffective for some scientific problems, such as high-dimensional partial differential equations (PDEs), inverse PDE problem, stochastic PDE problem, etc \cite{R14,R15}. Besides, most of these multi-scale methods require high-quality mesh discretization to computational domain and the numerical accuracy of these multi-scale methods has a strong dependence on grid quality. Hence, researchers further explore and develop effective machine learning (ML) or deep learning (DL) methods to solve partial differential equations (PDEs).

To the best of our knowledge, machine learning and deep learning methods have been successfully applied to computer vision, natural language processing, medical diagnosis and autonomous driving, among others. In 1998, Lagaris et al. \cite{R47} was the first to solve ordinary and partial differential equations by utilizing artificial neural networks (ANNs). However, owing to the limitations of computing resources and efficient algorithms in that era, this study could not attract much attention. In keeping with the recent and extraordinary success of ML methods, ML methods have also played an increasingly important role in solving PDEs. Specifically, Rudy et al. \cite{R14} proposed an efficient and robust sparse regression method, which can discover the governing PDEs of a given system via time-series measurements in the spatial domain. In reference \cite{R15}, a numerical framework using observation data and deep neural networks (DNN) are developed by Qin et al. for approximating unknown dynamical systems. Wang et al. \cite{R41} combined deep learning techniques and local multi-scale model reduction methodologies to design novel multi-layer neural networks to predict flow dynamics in porous media. In reference \cite{R16}, combining multi-layer perceptron, convolutional neural networks and long short-term memory structures, Kiyani et al. presented data-driven architectures to discover the nonlinear equations of motion for phase-field models with high-order derivatives. The above-mentioned ML methods are known as data-driven methods, whose accuracy depends heavily on the capability to generate large amounts of training data. In recent years, scientists devised some novel ML approaches to solve PDEs, whose core feature is to utilize the physical laws described by PDEs as neural network constraints. By using the automatic differentiation technique in neural networks, Raissi et al. \cite{R17} developed physics-informed neural networks (PINNs) and systematically elaborated on the core idea of this ML method. The core idea and mechanism of PINNs are to incorporate loss functions based on the residual errors of the governing equations, boundary conditions and initial conditions, and then minimize the residuals for training deep neural networks with the minimal or complete absence of data to obtain the solutions of PDEs. The numerical experiments showed the accuracy and efficiency of the PINNs method for solving several classical PDEs including forward and inverse problems. Attracting by the excellent computing performance of PINNs, plenty of improved and subsequent works have also appeared on this subject. Lu et al. \cite{R21} presented a residual-based adaptive refinement (RAR) trick to improve the efficiency of PINNs and a Python library called DeepXDE, which can be used to solve computational science and engineering problems. Yu et al. \cite{R22} proposed a gradient-enhanced physics-informed neural networks (gPINNs), which can improve the numerical accuracy of PINNS. In the meanwhile, its validity is verified by both forward and inverse PDE problems. In addition, Meng et al. \cite{R23} put forward a parareal physical information neural network (PPINN) for time-dependent PDEs. By decomposing a long-time problem into many independent short-time problems, the proposed PPINN approach may achieve a significant speed-up for long-time integration of PDEs. Jagtap et al. \cite{R24,R25} established a generalized space-time domain decomposition framework, named eXtended PINNs (XPINNs), to solve nonlinear partial differential equations (PDEs) on arbitrary complex-geometry domains. The proposed XPINNs with adaptive activation functions and dynamical weights techniques are successfully employed to solve the inverse supersonic flow problems involving expansion and compression waves with high discontinuities. Mao et al. \cite{R26} simulated high-speed aerodynamic flow by solving the forward and inverse Euler equations using PINNs. Moreover, applying a conceptual framework similar to the soft self-attention mechanism, McClenny et al. \cite{R27} developed a self-adaptive physics-informed neural networks (SA-PINNs). Numerical experiments on several linear and nonlinear PDE demonstrated that the proposed SA-PINNs are more accurate than other PINNs algorithms. In the field of solid mechanics, Henkes et al. \cite{R29} presented novel conservative PINNs (cPINNs), which have been successfully employed in solving nonlinear stress and displacement fields of heterogeneous materials with sharp phase interfaces. In addition to the above PINNs and their improved versions, scientists further developed other physics-driven ML approaches. E and Yu \cite{R18} first proposed a deep Ritz method (DRM) for solving variational problems of PDEs with residual neural network and applied it to simulate high-dimensional PDEs. In the same year, Sirignano et al. established a deep Galerkin method (DGM) for solving high-dimensional PDEs in reference \cite{R19}. The investigated numerical experiments demonstrated the exceptional approximation capability of neural networks for a class of quasi-linear parabolic PDEs. Furthermore, Samaniego et al. \cite{R20} proposed a deep energy method (DEM) and the application of this approach in computational mechanics was readily explored. In reference \cite{R28}, Peng et al. proposed a novel non-gradient deep learning method that can be used to solve high-dimensional elliptic equations such as Poisson equation, Fokker-Planck equation and advection-diffusion equation. In addition, Chen et al. \cite{R36} developed a deep mixed residual method (MIM) to compute the high-order PDEs via transforming a high-order PDE into a first-order system. In 2023, a deep double Ritz method (D$^2$RM) was provided by Uriarte et al. \cite{R48}, which combines two neural networks for approximating trial functions and optimal test functions along a nested double Ritz minimization strategy. Numerical experiments on different diffusion and convection problems support the robustness of the proposed D$^2$RM method. Yang et al. \cite{R51} employed deep Ritz method to investigate high-dimensional fractional-order differential equations with the right Riemann-Liouville fractional derivative, the left Caputo fractional derivative and the boundary value conditions. For the elliptic problems with singular sources, Lai et al. \cite{R52} developed a new shallow Ritz method using a shallow neural network with only one hidden layer, which significantly reduces the training costs. In reference \cite{R53}, Wang et al. proposed a mesh-free method to solve interface problems using the deep learning approach. The elliptic PDE with a discontinuous and high-contrast coefficient, and linear elasticity equation with discontinuous stress tensor were successfully simulated by the proposed ML method. Based on classical finite element error analysis, Minakowski et al. \cite{R54} derived the prior and posterior error estimates of deep Ritz method for Laplace and Stokes problems. In a summary, enormous ML and DL approaches has blossomed and widely employed in scientific computing and engineering application.

Although many neural network-based ML and DL approaches have been developed for solving mathematical problems, there are only a few work about multi-scale problems. Leung et al. \cite{R32} developed a neural homogenization-based PINN (NH-PINN), which combines the homogenization method and physics-informed neural networks to solve multi-scale problems. Xu and Li et al. \cite{R30,R31} proposed two types of multi-scale deep neural networks (MscaleDNNs) and ameliorated these MscaleDNNs via designing a new activation function. It is worth mentioning that these approaches have been applied to solve several complex elliptic equations. Zhang and Li et al. \cite{R33} combined subspace decomposition ideas and multi-scale DNN algorithms to solve multi-scale elliptic PDEs. Nevertheless, multi-scale problems considered in the above-mentioned studies are limited to the governing equations with continuous coefficients, which can not be directly employed to multi-scale PDEs of authentic composite materials with highly discontinuous coefficients. To fill this technique gap, by approximating the discontinuity coefficient using a continuous function, Chatzigeorgiou et al. \cite{R37} developed physically informed deep homogenization neural network (DHN) by integrating the periodic homogenization framework and physically informed neural networks (PINNs) to simulate the temperature field and temperature gradient distribution of unidirectional multiphase/multi-inclusion composites. But the presented DHN were only applied to compute composite materials with low-contrast material parameters and did not consider the numerical simulation of high-contrast composites. Furthermore, the foregoing multi-scale ML approaches are established based on asymptotic homogenization framework, which can not provide the adequate accuracy for engineering computation and capture the oscillatory behavior of heterogeneous materials at micro-scale. Additionally, classical neural network-based learning approaches have been theoretically and numerically demonstrated that they will fail in solving multi-scale problems due to so-called Frequency Principle (F-Principle) \cite{R32,R49,R50}. Also, direct ML and DL simulations of multi-scale problems will invest a huge amount of computing resources. Based on the above issues, it is of great practical value and theoretical significance to develop novel ML and DL approaches for multi-scale problems of authentic composite materials.

In this paper, we now propose a novel higher-order multi-scale deep Ritz method by seamlessly combining higher-order multi-scale method and deep Ritz method to solve multi-scale thermal transfer equation of authentic composite materials with highly oscillatory and discontinuous coefficients. As a demonstration, we implement HOMS-DRM to simulate the thermal transfer problems of composite materials with different types, including varying inclusion morphology of composites, high-contrast composites, and three-phase composites. The main contributions of this study include:
\begin{enumerate}
    \item A mesh-free and data-free higher-order multi-scale deep Ritz method (HOMS-DRM) is proposed, which tackles the intractable or even ineffective issues for solving multi-scale problems by classical neural network-based learning approaches.
    \item The proposed HOMS-DRM greatly improves the numerical accuracy for simulating multi-scale heat transfer equations with highly discontinuous coefficients of authentic and high-contrast composite materials. The microscopic oscillatory information of authentic and high-contrast composite materials can be precisely captured.
	\item Compared with classical neural network-based learning approaches, the presented HOMS-DRM can greatly economize the computing resources and improve computational efficiency when solving multi-scale problems.
	\item The theoretical convergence of the proposed HOMS-DRM method is rigorously demonstrated under appropriate assumptions.
\end{enumerate}

The remainder of this study is organized as follows: In Section \ref{sec:2}, the detailed establishment of higher-order multi-scale computational model is presented for steady-state thermal transfer problems of composite materials. After this, we review the basics of deep Ritz method and devise higher-order multi-scale deep Ritz method to accurately solve higher-order multi-scale computational model by precisely enforcing Dirichlet boundary conditions in Section \ref{sec:3}. Moreover, we also present its algorithm framework in detail. In Section \ref{sec:4}, we rigorously prove the convergence of HOMS-DRM under certain assumptions. Then some numerical examples are conducted in Section \ref{sec:5}, to verify the computational accuracy and efficiency of the presented HOMT-DRM. Finally, some meaningful conclusions and an outlook to necessary future work are summarized in Section \ref{sec:6}.

In this work, we apply Einstein summation convention to simplify repetitious indices.
\section{The higher-order multi-scale computational model for heat transfer equations in composite materials}
\label{sec:2}
\begin{figure}
	\centering
	\includegraphics[width=110mm]{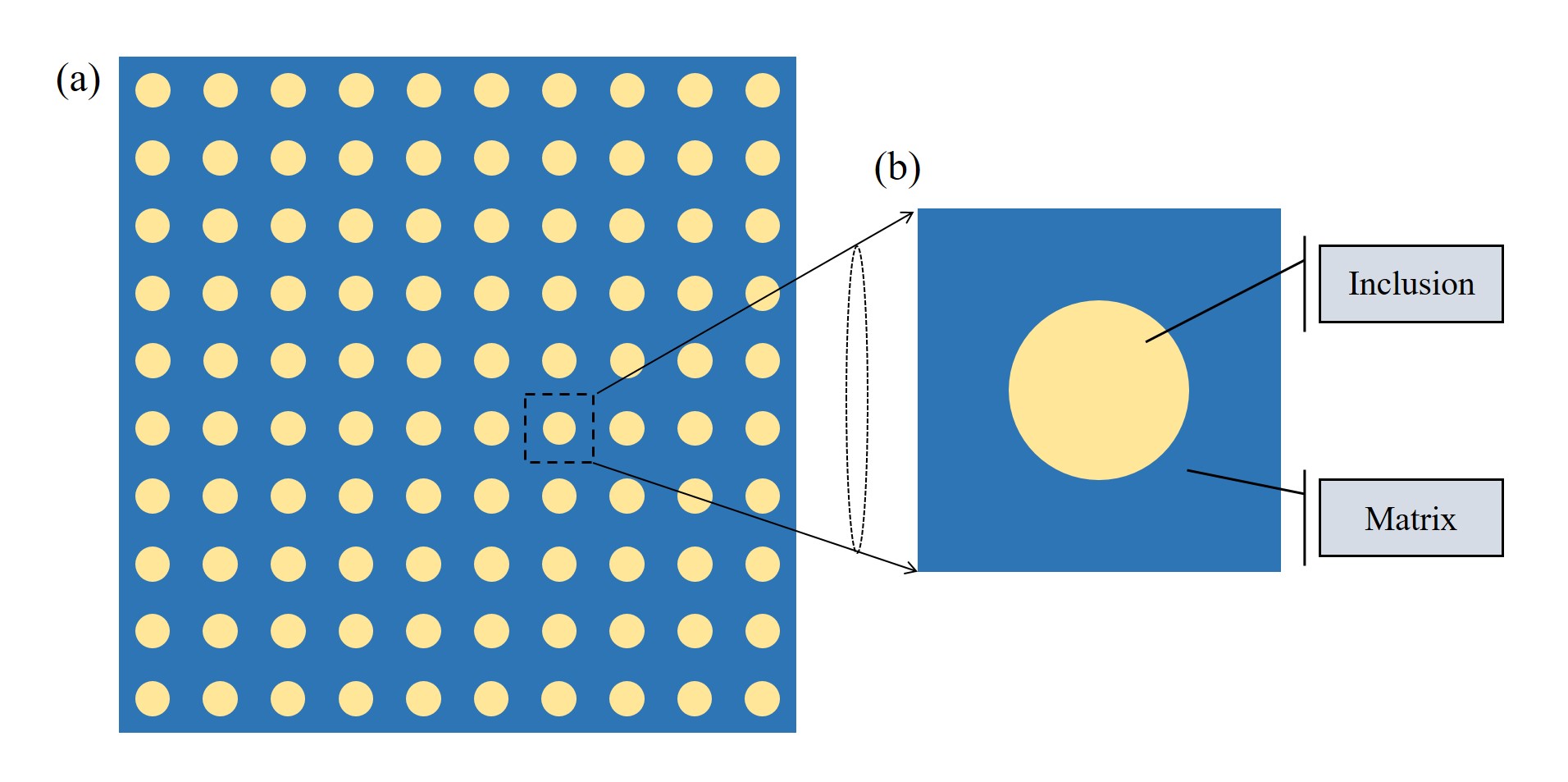}
	\caption{The diagram of periodic composite material. (a) The macroscopic structure $\Omega$; (b) The microscopic unit cell $Y$.}\label{fig:1}
\end{figure}

In this paper, we consider the steady-state thermal transfer problem of periodic composite materials. The investigated composite materials is schematically displayed in Fig.\hspace{1mm}\ref{fig:1}.

Based on the classical thermal transfer theory, the steady-state governing equation of periodic composite materials is presented as follows.
\begin{equation}
\label{eq:1}
\left\{  \begin{aligned}
&{ - \frac{\partial }{{\partial {x_i}}}\Big( {a_{ij}^\xi (\bm{x})\frac{{\partial {u^\xi }(\bm{x})}}{{\partial {x_j}}}} \Big) = f(\bm{x})},\quad{{\rm{in}}}\ \Omega,  \\
& {{u^\xi }(\bm{x}) = 0}, \quad {{\rm{on}}}\ {\partial \Omega }, \\
\end{aligned} \right.
\end{equation}
in which $\Omega $ is a bounded convex domain in ${R ^2}$  with an outer boundary $\partial \Omega $. The superscript $\xi $  indicates the characteristic size of microscopic unit cell $Y $; $a_{ij}^\xi (\bm{x})$ is the second-order thermal conductivity tensor; $f(\bm{x})$ is the internal thermal source; ${u^\xi}(\bm{x})$ is the unresolved temperature field of multi-scale heat conduction problem (1).

Next, we briefly present the higher-order multi-scale computational model for heat transfer equations in composite materials with periodic configurations. Firstly, we define variable $\bm{y}$ as the microscopic coordinate of microscopic unit cell $Y$ satisfying:
\begin{equation}
\label{eq:2}
\bm{y} = \frac{\bm{x}}{\xi} - \left[ {\frac{\bm{x}}{\xi}} \right] \in Y= {\left[ {0,1} \right]^2}.
\end{equation}
With the preceding definition (2), the thermal conductivity coefficient $a_{ij}^\xi (\bm{x})$ can be rewritten as ${a_{ij}}(\bm{y})$. In addition, we have the following chain rule connecting micro-scale and macro-scale.
\begin{equation}
\label{eq:3}
\frac{\partial \Phi^\xi(\bm{x})}{{\partial {x_i}}}=\frac{\partial \Phi(\bm{x},\bm{y})}{{\partial {x_i}}} + \frac{1}{\xi }\frac{\partial \Phi(\bm{x},\bm{y})}{{\partial {y_i}}}.
\end{equation}

Furthermore, classical multi-scale asymptotic analysis is implemented to multi-scale heat conduction problem (1) \cite{R42}. According to existing reference \cite{R42}, the higher-order multi-scale asymptotic solution $u^{( {2\xi })}(\bm{x})$ for multi-scale heat conduction problem (1) is derived as below.
\begin{equation}
\label{eq:4}
{u^{\left( {2\xi } \right)}}(\bm{x})={u_0}(\bm{x}) + \xi {N_{{\alpha _1}}}(\bm{y})\frac{{\partial {u_0}(\bm{x})}}{{\partial {x_{{\alpha _1}}}}} + {\xi ^2}{N_{{\alpha _1}{\alpha _2}}}(\bm{y})\frac{{{\partial ^2}{u_0}(\bm{x})}}{{\partial {x_{{\alpha _1}}}\partial {x_{{\alpha _2}}}}}.
\end{equation}
In the above formula (4), ${u_0}(\bm{x})$ is denoted as macroscopic homogenized solution for multi-scale heat conduction problem (1) defined on the macroscopic structure $\Omega $, ${N_{{\alpha _1}}}(\bm{y})$ and ${N_{{\alpha _1}{\alpha _2}}}(\bm{y})$ $({\alpha _1},{\alpha _2} = 1,2) $ are denoted as the lower-order (first-order) and higher-order (second-order) auxiliary cell functions defined on microscopic unit cell $Y$, respectively.

After that, substituting higher-order multi-scale solution $u^{( {2\xi })}(\bm{x})$ into multi-scale heat conduction problem (1) and equating the terms with the same power $\xi$ by chain rule (3), we can derive auxiliary cell problems for computing lower-order cell functions ${N_{{\alpha _1}}}(\bm{y})$ from the $O(\xi^{-1})$-order equation as follows.
\begin{equation}
\label{eq:5}
\left\{  \begin{aligned}
&{ - \frac{\partial }{{\partial {y_i}}}\Big( {{a_{ij}}(\bm{y})\frac{{\partial {N_{{\alpha _1}}}(\bm{y})}}{{\partial {y_j}}}} \Big) = \frac{{\partial {a_{i{\alpha _1}}}(\bm{y})}}{{\partial {y_i}}}}, \quad {{\rm{in}}}\ Y,  \\
& {{N_{{\alpha _1}}}(\bm{y}) = 0}, \quad {{\rm{on}}}\ {\partial Y}. \\
\end{aligned} \right.
\end{equation}
Afterwards, applying the homogenized operator to both sides of the $O(\xi^{0})$-order equation, we
shall derive the macroscopic homogenized equation as follows.
\begin{equation}
\label{eq:8}
\left\{  \begin{aligned}
&{ - \frac{\partial }{{\partial {x_i}}}\Big( {a_{ij}^0\frac{{\partial {u_0 }(\bm{x})}}{{\partial {x_j}}}} \Big) = f(\bm{x})},\quad{{\rm{in}}}\ \Omega,  \\
& {{u_0}(\bm{x}) = 0}, \quad {{\rm{on}}}\ {\partial \Omega }. \\
\end{aligned} \right.
\end{equation}
In macroscopic heat conduction equation (6), macroscopic homogenized coefficient $a_{ij}^0$ is computed by the following formula.
\begin{equation}
\label{eq:6}
a_{ij}^0 = \frac{1}{{|Y|}}\int_Y {\Big( {{a_{ij}}(\bm{y}) + {a_{ik}}(\bm{y})\frac{{\partial {N_j}(\bm{y})}}{{\partial {y_k}}}} \Big)} dY.
\end{equation}
Finally, the following auxiliary cell problems for uniquely determining higher-order cell functions ${N_{{\alpha _1}{\alpha _2}}}(\bm{y})$ are established from the $O(\xi^{0})$-order equation.
\begin{equation}
\label{eq:7}
\left\{  \begin{aligned}
& - \frac{\partial }{{\partial {y_i}}}\Big( {{a_{ij}}(\bm{y})\frac{{\partial {N_{{\alpha _1}{\alpha _2}}}(\bm{y})}}{{\partial {y_j}}}} \Big) = \frac{\partial }{{\partial {y_i}}}\left( {{a_{i{\alpha _2}}}(\bm{y}){N_{{\alpha _1}}}(\bm{y})} \right) \\
&+ {a_{{\alpha _2}j}}(\bm{y})\frac{{\partial {N_{{\alpha _1}}}(\bm{y})}}{{\partial {y_j}}} + {a_{{\alpha _1}{\alpha _2}}}(\bm{y}) - a_{{\alpha _1}{\alpha _2}}^0{\rm{,     \quad in }}\ Y,\\
&	{{N_{{\alpha _1}{\alpha _2}}}(\bm{y}) = 0,\quad {\rm{ on }}\ \partial Y}.
\end{aligned} \right.
\end{equation}

Summing up, the higher-order multi-scale computational model for heat transfer equations in composite materials is established by introducing higher-order cell functions $N_{{\alpha _1}{\alpha _2}}(\bm{y})$ and corresponding higher-order corrected term $\displaystyle N_{{\alpha _1}{\alpha _2}}(\bm{y})\frac{{\partial ^2}{u_0}(\bm{x})}{\partial x_{{\alpha _1}}\partial x_{{\alpha _2}}}$, which exhibiting not only the exceptional numerical accuracy, but also the less computational expense. In addition, it is worth emphasizing that we replace classical periodic boundary condition for microscopic cell functions with homogeneous Dirichlet boundary condition \cite{R42,R43,R44}. For microscopic unit cell with geometric symmetry, it can be proved that Dirichlet boundary condition is equivalent to classical periodic boundary condition \cite{R42,R43,R44}.
\section{Higher-order multi-scale deep Ritz method (HOMS-DRM)}
\label{sec:3}
In this section, we introduce our proposed approach: higher-order multi-scale deep Ritz method (HOMS-DRM) in detail. Based on previous researches \cite{R32,R49,R50}, we already have a clear understanding that classical neural network-based learning approaches fail in solving multi-scale problems with a satisfying accuracy due to Frequency Principle. Moreover, direct neural network-based simulations of multi-scale problems will consume a tremendous amount of computing resources. Furthermore, for authentic composite materials, their thermal conductivity in matrix and inclusion components are strongly discontinuous and piecewise constants, which is not derivable at the interface layer of authentic composite materials. The deep Ritz method possesses natural superiority to resolve this deficiency by transforming derivative to integral via the
principle of minimum potential energy. The aforementioned challenges urgently motivate us to develop novel computational approach for multi-scale problems, which can inherit the respective superiority of higher-order multi-scale method and deep Ritz method.

The core idea of our proposed HOMS-DRM consists of four main steps: (1) Establish higher-order multi-scale computational model for original multi-scale problems of authentic composite materials and decompose the hard multi-scale problems into easier and solvable macroscopic homogenized problems and microscopic cell problems without multi-scale property. (2) Solve the microscopic cell problems by improved deep Ritz method and evaluate macroscopic homogenized coefficients of authentic composite materials. (3) Solve macroscopic homogenized equation by using improved deep Ritz method. (4) Assemble higher-order multi-scale solutions for multi-scale problems of authentic composite materials. Now, we start to present higher-order multi-scale deep Ritz method in detail.
\subsection{The basic theory of deep Ritz method}
Briefly, we will first introduce some basic concepts of deep Ritz method (DRM), which is employed as neural network-
based technique for solving PDEs without multi-scale property. In this study, we focus on applying hard-constraint boundary conditions to deep Ritz method since macroscopic homogenized equation and microscopic cell problems of higher-order multi-scale computational model all are imposed with Dirichlet boundary condition. In principle, DRM is a neural network-based approach that uses artificial neural networks to approximate the solutions of PDEs by minimizing the potential energy of the investigated PDEs.

In this study, solving lower-order cell problems (5) are employed as an example to introduce the principal idea of DRM in detail. Based on the theoretical framework of DRM, we first assume ${\hat N_{{\alpha _1}}}({\bm{y}})$ is the DRM solution of lower-order cell functions ${N_{{\alpha _1}}}(\bm{y})$, which is approximated by using a deep residual network $F( {{\bm{y}};{\bm{\theta}}})$ \cite{R18}. The deep residual network $F( {{\bm{y}};{\bm{\theta}}})$ own input variable $\bm{y}=\left(y_1,y_2\right)$ and a set of all network parameters $\bm{\theta}\in R^D$ including weight and bias parameters, where $D$ is the total number of network parameters. To the best of our knowledge, the residual neural networks were first introduced in large-scale visual recognition where they exhibited improved performance that traditional neural networks \cite{R35}. In this work, the employed deep residual network consists of full connected layer and the residual block. The mathematical representation of the full connected layer can be written as below.
\begin{equation}
{{\bm{z}}^o} = \sigma \left( {\bm{W} \cdot {{\bm{z}}^i} + {\bm{b}}} \right),
\end{equation}
in which $\bm{W}$ and $\bm{b}$ respectively refer to the weight and bias parameters of the fully connected layer. ${\bm{z}}^i$ and ${\bm{z}}^o$ are the input and output of the layer, respectively. $\sigma$ is well-known as nonlinear activation function. The commonly employed activation functions include ReLU function: $\sigma \left( z \right) = \max \left( {0,z} \right)$, Sigmoid function: $\sigma \left( z \right) = 1/\left( {1 + {e^{ - z}}} \right)$ and Tanh function: $\sigma \left( z \right) = \left( {{e^z} - {e^{ - z}}} \right)/\left( {{e^z} + {e^{ - z}}} \right)$. Additionally, for some specific problems, researchers have developed some new activation functions, such as $ReLU^3$ function: $\sigma \left( z \right) = \max \left( {0,z^3} \right)$ \cite{R34} and $sReLU$ function: $\sigma \left( z \right) = {\left( x \right)_ + }{\left( {1 - x} \right)_ + }$ \cite{R30}, etc. On the other hand, the key feature of residual blocks is their ability to capture and leverage residual information within the network architecture. In this work, each residual block comprises two fully connected layers. And the mathematical representation of the residual block is given as follows.
\begin{equation}
{{\bm{z}}^o} = {{\bm{z}}^i} +\sigma \left( {{\bm{W}^2} \cdot \left( {\sigma \left( {{\bm{W}^1} \cdot {{\bm{z}}^i} + {{\bm{b}}^1}} \right)} \right) + {{\bm{b}}^2}} \right).
\end{equation}

\begin{figure}
	\centering
	\includegraphics[width=140mm]{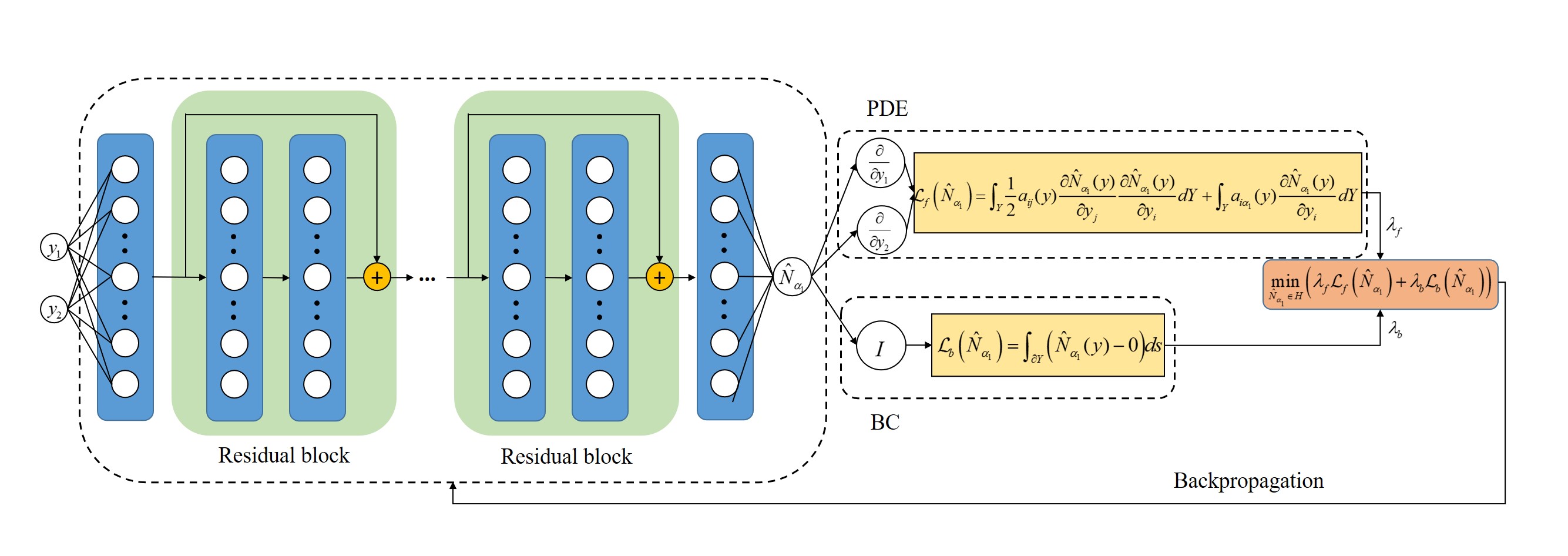}
	\caption{The schematic of DRM for solving lower-order cell functions ${N_{{\alpha _1}}}(\bm{y})$.}\label{fig:2}
\end{figure}

Next, the loss function of deep residual network $F( {{\bm{y}};{\bm{\theta}}})$ is defined by using the outputs of the network and their derivatives on the basis of the principle of minimum potential energy. As a result, the loss function of deep residual network $F( {{\bm{y}};{\bm{\theta}}})$ is presented as follows.
\begin{equation}
{{\cal L}_F} = {\lambda _f}{{\cal L}_f} + {\lambda _b}{{\cal L}_b},
\end{equation}
where ${{\cal L}_f}$ and ${{\cal L}_b}$ represent the losses of solving PDEs coming from inner domain and boundary respectively, and are presented as follows.
\begin{equation}
\begin{aligned}
&{{\cal L}_f} = \int_Y \frac{1}{2}{{a_{ij}}} (\bm{y})\frac{{\partial {{\hat N}_{{\alpha _1}}}(\bm{y})}}{{\partial {y_j}}}\frac{{\partial {{\hat N}_{{\alpha _1}}}(\bm{y})}}{{\partial {y_i}}}dY + \int_Y {{a_{i{\alpha _1}}}} (\bm{y})\frac{{\partial {{\hat N}_{{\alpha _1}}}(\bm{y})}}{{\partial {y_i}}}dY\\
&\approx \frac{1}{{{n_{in}}}}\sum\limits_{k = 1}^{{n_{in}}} {\Big[ {\frac{1}{2}\big( {{a_{ij}}({\bm{y}_k})\frac{{\partial {{\hat N}_{{\alpha _1}}}({\bm{y}_k})}}{{\partial {y_j}}}\frac{{\partial {{\hat N}_{{\alpha _1}}}({\bm{y}_k})}}{{\partial {y_i}}}} \big) + \big( {{a_{i{\alpha _1}}}({\bm{y}_k})\frac{{\partial {{\hat N}_{{\alpha _1}}}({\bm{y}_k})}}{{\partial {y_i}}}} \big)} \Big]},
\end{aligned}
\end{equation}
\begin{equation}
{{\cal L}_b} = \int_{\partial Y} {\big( {{{\hat N}_{{\alpha _1}}}(\bm{y}) - 0} \big)} ds \approx \frac{1}{{{n_{bd}}}}\sum\limits_{k = 1}^{{n_{bd}}} {\big( {{{\hat N}_{{\alpha _1}}}({\bm{y}_k}) - 0} \big)}.
\end{equation}
Moreover, ${\lambda _f}$ and ${\lambda _b}$ in loss function (9) denote the penalty or weight coefficients, which can enable deep residual network to better approximate the boundary term of PDEs. Furthermore, in the above equations (10) and (11), ${n_{in}}$ and ${n_{bd}}$ are defined as the discretization number of interior points and boundary points, respectively. In a summary, the schematic diagram of deep Ritz method based on residual neural networks is shown in Fig.\hspace{1mm}\ref{fig:2}.

To handle the Dirichlet boundary of macroscopic homogenized problems and microscopic cell problems, another strategy has been proposed by directly imposing the hard-constraint boundary condition into neural networks \cite{R29,R55}. Taking lower-order cell functions ${N_{{\alpha _1}}}(\bm{y})$ as example, Dirichlet boundary condition ${N_{{\alpha _1}}}(\bm{y})|_{\bm{y}
\in\partial Y}=0$ can be coercively imposed by modifying original deep residual network $F( {{\bm{y}};{\bm{\theta}}})$ as new deep neural network $G\left( \bm{y};\bm{\theta} \right) = {y_1} \cdot {y_2} \cdot \left( {1 - {y_1}} \right) \cdot \left( {1 - {y_2}} \right) \cdot F\left( \bm{y};\bm{\theta} \right)$. Furthermore, it should be emphasized that modified neural network $G\left( \bm{y};\bm{\theta} \right)$ only considers the inner loss for solving PDEs if hard-constrained boundary conditions are applied. Eventually, neural network-based approximate solution ${\hat N_{{\alpha _1}}}({\bm{y}})$ is obtained by minimizing new loss function of modified neural network $G\left( \bm{y};\bm{\theta} \right)$ as below.
\begin{equation}
{\bm{\theta} ^*} = \mathop {{\mathop{\rm argmin}\nolimits} }\limits_{{\bm{\theta }} \in {R ^D}} {{\cal L}_G}(\bm{y};\bm{\theta}).
\end{equation}
\subsection{The theoretical framework of higher-order multi-scale deep Ritz method}
\begin{figure}
	\centering
	\includegraphics[width=140mm]{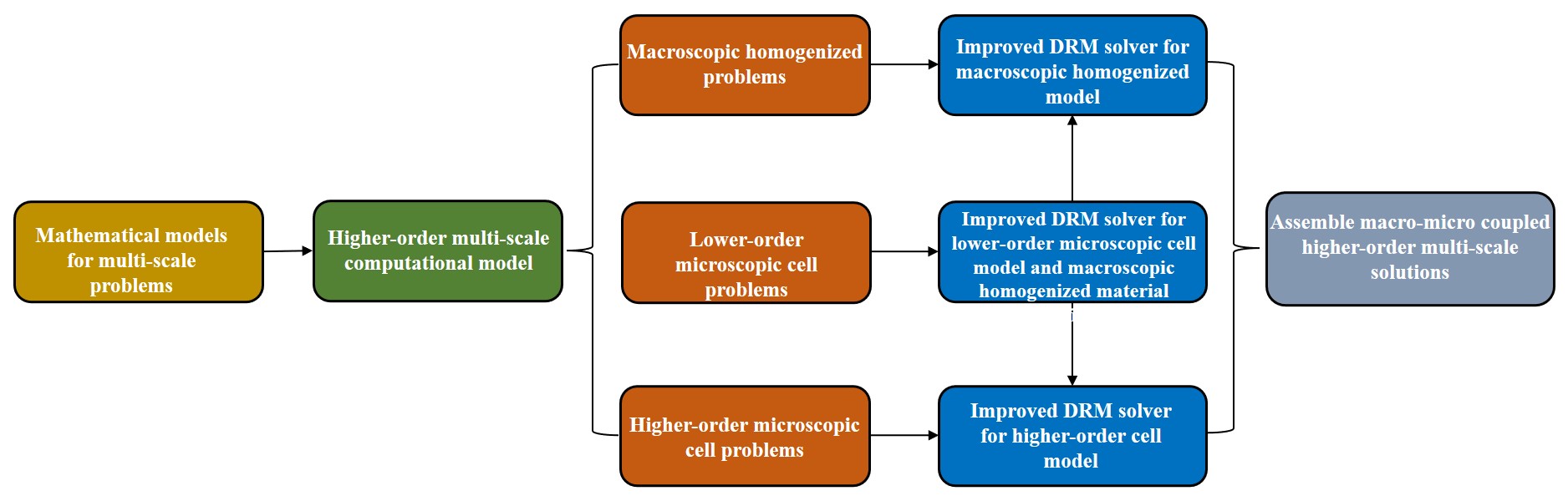}
	\caption{The computational framework of HOMS-DRM.}\label{fig:3}
\end{figure}

As mentioned at the beginning of this section, it is extremely difficult to directly solve multi-scale problems of authentic composite materials via classical neural network-based learning approaches. Thereby, we first establish higher-order multi-scale computational model for original multi-scale problems of authentic composite materials. The presented higher-order multi-scale computational model decompose original multi-scale problems into macroscopic homogenized problems and microscopic cell problems without multi-scale property. Next, improved deep Ritz method with hard-constraint boundary condition are developed for mesh-free solving macroscopic homogenized problems and microscopic cell problems of authentic composite materials, respectively. Finally, high-accuracy multi-scale solutions for heat transfer equations in authentic composite materials are obtained based on higher-order multi-scale computational model and automatic differentiation technique, which can accurately capture the microscopic oscillatory information of authentic composite materials. The detailed computational framework of HOMS-DRM is displayed in Fig.\hspace{1mm}\ref{fig:3}.

Furthermore, we evaluate the macroscopic homogenized coefficients $a_{ij}^0$ by formula (7), which demands to calculate the derivatives of lower-order microscopic cell functions and do numerical integration. Owing to the obtained deep learning solutions of lower-order microscopic cell functions and automatic differentiation package in modern deep learning software, we can give an accurate approximation of macroscopic homogenized coefficient in a relatively easier way and denote the resulting macroscopic homogenized coefficient as $\hat a_{ij}^0$. After that, we successively establish loss functions for macroscopic homogenized problems (6) and higher-order microscopic cell problems (8) by performing the principle of minimum potential energy as follows.
\begin{equation}
\label{eq:11}
{{\cal L}_f}=\int_\Omega  \dfrac{1}{2}{\hat a_{ij}^0} \frac{{\partial {\hat u_0}(\bm{x})}}{{\partial {x_j}}}\frac{{\partial {\hat u_0}(\bm{x})}}{{\partial {x_i}}}d\Omega-\int_\Omega  f (\bm{x}){\hat u_0}(\bm{x})d\Omega.
\end{equation}
\begin{equation}
{{\cal L}_b} = \int_{\partial \Omega} {\big( {{{\hat u}_0}(\bm{x}) - 0} \big)} ds.
\end{equation}
\begin{equation}
\begin{aligned}
&{{\cal L}_f}=\int_Y \dfrac{1}{2}{{a_{ij}}} (\bm{y})\frac{{\partial {\hat N_{{\alpha _1}{\alpha _2}}}(\bm{y})}}{{\partial {y_j}}}\frac{{\partial {\hat N_{{\alpha _1}{\alpha _2}}}(\bm{y}) }}{{\partial {y_i}}}dY\\
&-{\rm{ }}\int_Y {\big( {{a_{{\alpha _2}j}}(\bm{y})\frac{{\partial {\hat N_{{\alpha _1}}}(\bm{y})}}{{\partial {y_j}}} + {a_{{\alpha _1}{\alpha _2}}}(\bm{y}) - \hat a_{{\alpha _1}{\alpha _2}}^0} \big)} {\hat N_{{\alpha _1}{\alpha _2}}}(\bm{y}) dY \\
&+\int_Y {{a_{i{\alpha _2}}}} (\bm{y}){N_{{\alpha _1}}}(\bm{y})\frac{{\partial {\hat N_{{\alpha _1}{\alpha _2}}}(\bm{y}) }}{{\partial {y_i}}}dY.
\end{aligned}
\end{equation}
\begin{equation}
{{\cal L}_b} = \int_{\partial Y} {\big( {{{\hat N}_{{\alpha _1}{\alpha _2}}}(\bm{y}) - 0} \big)} ds.
\end{equation}
In above-mentioned loss functions (15)-(18), ${\hat u_0}(\bm{x})$ and ${\hat N_{{\alpha _1\alpha _2}}}({\bm{y}})$ are the DRM solutions of macroscopic homogenized solutions ${u_0}(\bm{x})$ and higher-order microscopic cell functions ${N_{{\alpha _1\alpha _2}}}(\bm{y})$, which all are approximated by employing deep residual networks. Then, by utilizing improved deep Ritz method with hard-constraint boundary condition, we can naturally obtain high-accuracy numerical solutions for ${u_0}(\bm{x})$ and ${N_{{\alpha _1\alpha _2}}}(\bm{y})$. Lastly, higher-order multi-scale deep Ritz solution $\hat u^{\left( {2\xi } \right)}(\bm{x})$ for thermal transfer equation (1) of composite materials is obtained as follows.
\begin{equation}
{\hat u^{\left( {2\xi } \right)}}(\bm{x})={\hat u_0}(\bm{x}) + \xi {\hat N_{{\alpha _1}}}(\bm{y})\frac{{\partial {\hat u_0}(\bm{x})}}{{\partial {x_{{\alpha _1}}}}} + {\xi ^2}{\hat N_{{\alpha _1}{\alpha _2}}}(\bm{y})\frac{{{\partial ^2}{\hat u_0}(\bm{x})}}{{\partial {x_{{\alpha _1}}}\partial {x_{{\alpha _2}}}}}.
\end{equation}
\subsection{Algorithm procedures of higher-order multi-scale deep Ritz method}
In this subsection, we present the algorithm procedures of high-order multi-scale deep Ritz method (HOMS-DRM) in detail. The algorithm procedures consist of the following five primary steps.
\begin{enumerate}[Step 1:]
	\item Establish higher-order multi-scale computational model for multi-scale heat transfer problems and derive the detailed mathematical models for macroscopic homogenized solution ${u_0}(\bm{x})$, lower-order microscopic cell functions ${N_{\alpha_1}}(\bm{y})$ and higher-order microscopic cell functions ${N_{{\alpha_1}{\alpha_2}}}(\bm{y})$.
	\item Solve lower-order microscopic cell functions ${N_{\alpha_1}}(\bm{y})$ by using improved DRM with loss functions (12) and (13) to auxiliary cell problems (5) and calculate macroscopic homogenized material coefficients $a_{ij}^0$ by formula (7).
	\item Solve macroscopic homogenized solution ${u_0}(\bm{x})$ by applying improved DRM with loss functions (15) and (16) to macroscopic homogenized heat transfer problems (6).
	\item Solve higher-order microscopic cell functions ${N_{{\alpha_1}{\alpha_2}}}(\bm{y})$ by employing improved DRM with loss functions (17) and (18) to auxiliary cell problems (8).
	\item Employ formula (19) to assemble high-accuracy higher-order multi-scale asymptotic solution for multi-scale heat transfer problems of authentic composite materials. The spatial derivatives $\displaystyle\frac{{\partial {\hat u_0}(\bm{x})}}{{\partial {x_{{\alpha _1}}}}}$ and $\displaystyle\frac{{{\partial ^2}{\hat u_0}(\bm{x})}}{{\partial {x_{{\alpha _1}}}\partial {x_{{\alpha _2}}}}}$ in macro-micro coupled formula (19) are evaluated by automatic differentiation technique.
\end{enumerate}
\section{The proof of main convergence}
\label{sec:4}
In this section, we will prove the convergence of the proposed HOMS-DRM under the appropriate assumptions. There are four error sources in our presented HOMS-DRM. Firstly, the approximation error arises from the higher-order multi-scale model. Secondly, other two approximation errors originate from the improved DRM solver accuracy when solving microscopic cell problems and macroscopic homogenized equation, respectively. The another error comes from optimization error of determining neural network, however, this error is small compared with approximation errors. Hence, the proof of main convergence will only consider approximation errors and we shall exhibit in the numerical section that the computational error of our proposed method can be controlled.
\begin{lem}
According to the existing results of error estimation in \cite{R42,R44}, the following error estimation has been obtained between exact solution $u^{\xi}(\bm{x})$ and higher-order multi-scale solution $u^{(2\xi)}(\bm{x})$.
\begin{equation}
\Big\|u^{\xi}(\bm{x})-u^{(2\xi)}(\bm{x})\Big\|_{H^1(\Omega)}\leq C\xi,
\end{equation}
where $C$ is a positive constant irrespective of $\xi$, but dependent of $\Omega$.
\end{lem}
\begin{lem}
Based on the existing results of error estimation in \cite{R56}, there exists a network $U\in\mathcal{F}(\Phi)$ where $\Phi$ is the common activation function of the hidden units, which is uniformly 2-dense on compacts of $H^2(\bar\Omega)$. It means that for $u\in H^2(\bar\Omega)$ and $\forall\varepsilon>0$, it follows that
\begin{equation}
\mathop{\max}\limits_{a\le2}\mathop{\sup}\limits_{\bm{x}\in\Omega}|\partial _{\bm{x}}^au(\bm{x})-\partial _{\bm{x}}^aU(\bm{x};\bm{\beta})|<\varepsilon.
\end{equation}
\end{lem}
\begin{thm}
Under the above hypotheses lemmas 1 and 2, the following global error estimation are obtained.
\begin{equation}
{\begin{aligned}
\Big\|u^{\xi}(\bm{x})-\hat u^{(2\xi)}(\bm{x})\Big\|_{L^2(\Omega)}\leq C\xi+C\varepsilon+C\xi\varepsilon+C\xi^2\varepsilon.
\end{aligned}}
\end{equation}
\end{thm}
where $C$ is a positive constant irrespective of $\xi$, but dependent of $\Omega$.\\
$\mathbf{Proof:}$ Employing the triangle inequality and lemma 1, there firstly holds the following inequality.
\begin{equation}
{\begin{aligned}
&\quad\;\Big\|u^{\xi}(\bm{x})-\hat u^{(2\xi)}(\bm{x})\Big\|_{L^2(\Omega)}\\
&\leq \Big\|u^{\xi}(\bm{x})-u^{(2\xi)}(\bm{x})\Big\|_{L^2(\Omega)}+\Big\|u^{(2\xi)}(\bm{x})-\hat u^{(2\xi)}(\bm{x})\Big\|_{L^2(\Omega)}\\
&\leq \Big\|u^{\xi}(\bm{x})-u^{(2\xi)}(\bm{x})\Big\|_{H^1(\Omega)}+\Big\|u^{(2\xi)}(\bm{x})-\hat u^{(2\xi)}(\bm{x})\Big\|_{L^2(\Omega)}\\
&\leq C\xi+\Big\|u^{(0)}(\bm{x})-\hat u^{(0)}(\bm{x})\Big\|_{L^2(\Omega)}\\
&+\Big\|\xi{N_{\alpha_1}}({\bm{y}})\frac{\partial u^{(0)}({\bm{x}})}{\partial x_{\alpha_1}}-\xi{\hat N_{\alpha_1}}({\bm{y}})\frac{\partial \hat u^{(0)}({\bm{x}})}{\partial x_{\alpha_1}}\Big\|_{L^2(\Omega)}\\
&+\Big\|\xi^2{N_{\alpha_1\alpha_2}}({\bm{y}})\frac{\partial^2 u^{(0)}({\bm{x}})}{\partial x_{\alpha_1\alpha_2}}-\xi^2{\hat N_{\alpha_1\alpha_2}}({\bm{y}})\frac{\partial^2 \hat u^{(0)}({\bm{x}})}{\partial x_{\alpha_1\alpha_2}}\Big\|_{L^2(\Omega)}
\end{aligned}}
\end{equation}
After that, in virtue of the triangle inequality and lemma 2, we can derive three preliminary estimations of right sides of (23) as below.
\begin{equation}
{\begin{aligned}
\Big\|u^{(0)}(\bm{x})-\hat u^{(0)}(\bm{x})\Big\|_{L^2(\Omega)}\leq C\varepsilon.
\end{aligned}}
\end{equation}
\begin{equation}
{\begin{aligned}
&\quad\;\Big\|\xi{N_{\alpha_1}}({\bm{y}})\frac{\partial u^{(0)}({\bm{x}})}{\partial x_{\alpha_1}}-\xi{\hat N_{\alpha_1}}({\bm{y}})\frac{\partial \hat u^{(0)}({\bm{x}})}{\partial x_{\alpha_1}}\Big\|_{L^2(\Omega)}\\
&\leq \Big\|\xi\big[{N_{\alpha_1}}({\bm{y}})-{\hat N_{\alpha_1}}({\bm{y}})\big]\frac{\partial u^{(0)}({\bm{x}})}{\partial x_{\alpha_1}}\Big\|_{L^2(\Omega)}+\Big\|\xi{\hat N_{\alpha_1}}({\bm{y}})\big[\frac{\partial u^{(0)}({\bm{x}})}{\partial x_{\alpha_1}}-\frac{\partial \hat u^{(0)}({\bm{x}})}{\partial x_{\alpha_1}}\big]\Big\|_{L^2(\Omega)}\\
&\leq C\xi\varepsilon+C\xi\varepsilon\leq C\xi\varepsilon.
\end{aligned}}
\end{equation}
\begin{equation}
{\begin{aligned}
&\quad\;\Big\|\xi^2{N_{\alpha_1\alpha_2}}({\bm{y}})\frac{\partial^2 u^{(0)}({\bm{x}})}{\partial x_{\alpha_1\alpha_2}}-\xi^2{\hat N_{\alpha_1\alpha_2}}({\bm{y}})\frac{\partial^2 \hat u^{(0)}({\bm{x}})}{\partial x_{\alpha_1\alpha_2}}\Big\|_{L^2(\Omega)}\\
&\leq \Big\|\xi^2\big[{N_{\alpha_1\alpha_2}}({\bm{y}})-{\hat N_{\alpha_1\alpha_2}}({\bm{y}})\big]\frac{\partial^2 u^{(0)}({\bm{x}})}{\partial x_{\alpha_1\alpha_2}}\Big\|_{L^2(\Omega)}+\Big\|\xi^2{\hat N_{\alpha_1\alpha_2}}({\bm{y}})\big[\frac{\partial^2 u^{(0)}({\bm{x}})}{\partial x_{\alpha_1\alpha_2}}-\frac{\partial^2 \hat u^{(0)}({\bm{x}})}{\partial x_{\alpha_1\alpha_2}}\big]\Big\|_{L^2(\Omega)}\\
&\leq C\xi^2\varepsilon+C\xi^2\varepsilon\leq C\xi^2\varepsilon.
\end{aligned}}
\end{equation} 
Finally, substituting the above equalities (24)-(26) into the inequality (23), we hereby obtain the strongly convergent estimation as mentioned in (22), which guarantee the computational effectiveness of the proposed HOMS-DRM.
\section{Numerical examples and discussions}
\label{sec:5}
In this section, the proposed HOMS-DRM is employed to simulate the temperature fields of different types of composite materials in order to validate the computational capabilities of higher-order multi-scale deep Ritz method and corresponding numerical algorithm developed in this study. For the sake of simplicity, both matrix and inclusions herein are assumed to be isotropic. Moreover, we quantify the numerical accuracy of HOMS-DRM by introducing some relative error notations, which are defined as follows.
\begin{displaymath}
erro{r_1} = \frac{{{{\left\| {{u^D} - {u^\xi }} \right\|}_{{L_2}}}}}{{{{\left\| {{u^\xi }} \right\|}_{{L_2}}}}},erro{r_2} = \frac{{{{\left\| {{u^D} - {u^{(2\xi) }}} \right\|}_{{L_2}}}}}{{{{\left\| {{u^{(2\xi) }}} \right\|}_{{L_2}}}}},erro{r_3} = \frac{{{{\left\| {{u^{(2\xi) }} - {u^\xi }} \right\|}_{{L_2}}}}}{{{{\left\| {{u^\xi }} \right\|}_{{L_2}}}}},
\end{displaymath}
where ${u^D}$ is the numerical solution of the proposed HOMS-DRM, ${u^\xi }$ is the reference solution of the direct numerical simulation of fine finite element (DNS-FEM) and ${u^{(2\xi) }}$ is the reference solution of classical higher-order multi-scale (HOMS) method. For all numerical examples in this study, we employ the Tanh activation function and Adam optimizer to deep learning. Furthermore, all numerical experiments are executed on a desktop workstation equipped with an NVIDIA Quadro P620 graphics card, 16GB of internal storage, and a 2.60 GHz Intel Core i7-10750H CPU.
\subsection{Example 1: Thermal transfer simulation of composite materials with different thermal conductivities}
\label{subsec:1}
In this example, the proposed HOMS-DRM are employed to implement thermal transfer simulation of the composites with different thermal conductivities. The investigated composites structure is shown in Fig.\hspace{1mm}\ref{fig:1}, where $\Omega  = \left( {{x_1},{x_2}} \right) = \left[ {0,1} \right] \times \left[ {0,1} \right]$, $\xi  = 1/10$ and $Y = \left( {{y_1},{y_2}} \right) = \left[ {0,1} \right] \times \left[ {0,1} \right]$. In microscopic unit cell $Y$, the inclusion is a circle with central coordinate $\left( {0.5,0.5} \right)$ and radius $r = 0.3$. Furthermore, the internal thermal source $f(\bm{x})$ in this experiment is set as $10$ for thermal transfer simulation of composite materials. Then, the thermal transfer coefficients of matrix material are set as fixed ${a_{11}} = a_{22}=1$ and the thermal transfer coefficients in inclusion material are defined as varying 0.1, 0.01 and 0.001 respectively.

To implement thermal transfer simulation by our HOMS-DRM, we utilize a deep residual network with network structure 2-10-[10-10]$\times$4-10-1 to solve microscopic cell problems and macroscopic homogenized problem. The notation [10-10] represents one residual block and [10-10]$\times$n represents n residual blocks linked sequentially. Then randomly selecting 10,000 discrete points based on a uniform distribution within microscopic unit cell $Y$, we train the deep residual network for 30,000 epochs with Adam gradient descent (initial learning rate 0.005 and dynamic decrease to 0.95 times of previous one after every 500 epochs) to obtain deep learning solutions for lower-order microscopic cell functions $N_{\alpha_1}(\bm{y})$. After that, macroscopic homogenized coefficients are calculated by using 40,000 uniform-spaced points within lower-order microscopic cell functions $N_{\alpha_1}(\bm{y})$. Additionally, we continue to solve higher-order microscopic cell functions $N_{\alpha_1\alpha_2}(\bm{y})$ and macroscopic homogenized solution ${u_0}(\bm{x})$ with 10,000 random points. The deep learning solutions for higher-order microscopic cell functions $N_{\alpha_1\alpha_2}(\bm{y})$ and macroscopic homogenized solution ${u_0}(\bm{x})$ are obtained after 30,000 training epochs.

\begin{figure}
	\centering
	\includegraphics[width=130mm]{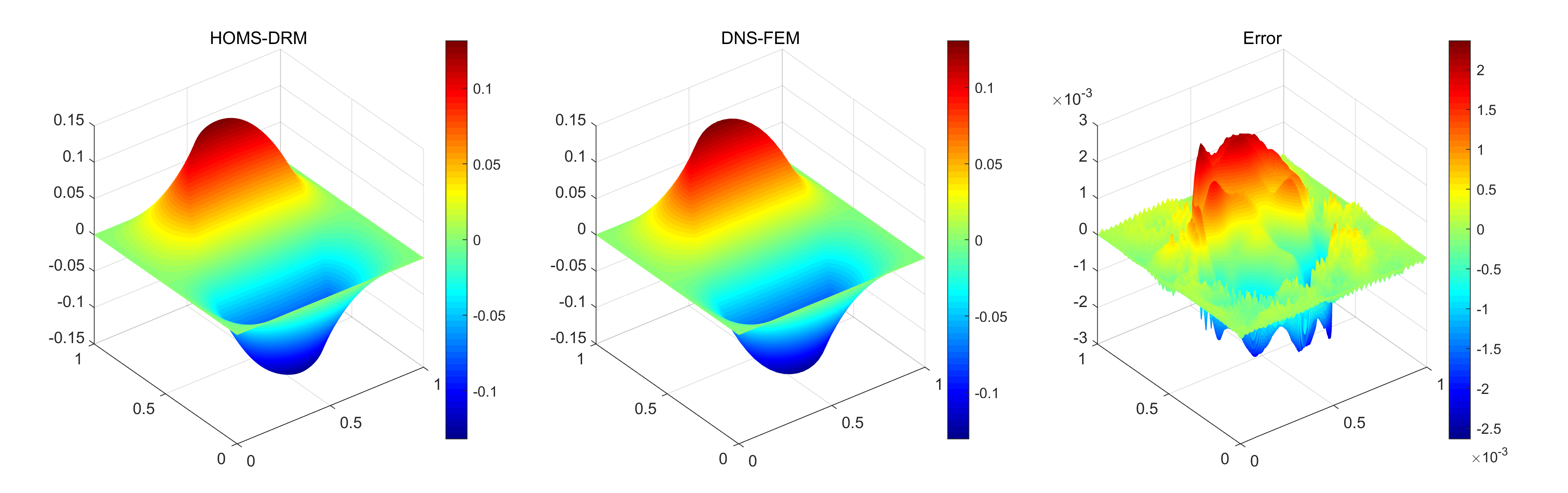}
	\caption{The computational results of lower-order microscopic cell function $N_1(\bm{y})$.}\label{fig:4}
\end{figure}

\begin{figure}
	\centering
	\includegraphics[width=130mm]{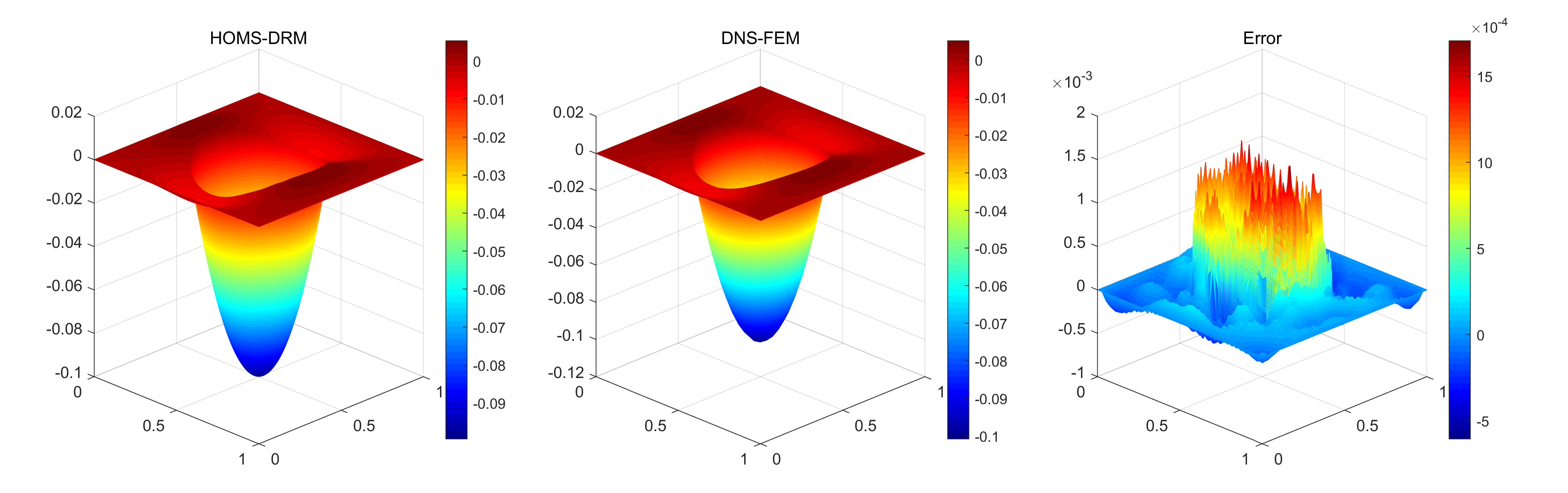}
	\caption{The computational results of higher-order microscopic cell function $N_{11}(\bm{y})$.}\label{fig:5}
\end{figure}

\begin{figure}
	\centering
	\includegraphics[width=130mm]{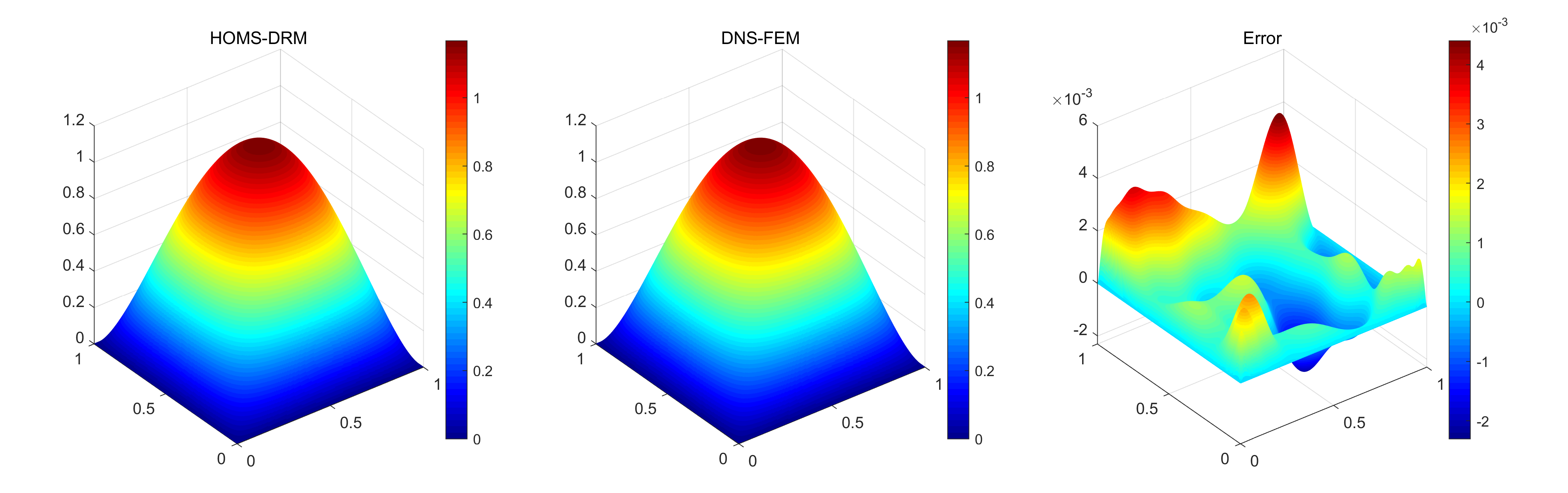}
	\caption{The computational results of macroscopic homogenized solution $u_0(\bm{x})$.}\label{fig:6}
\end{figure}

After numerical simulation, we present the comparisons of lower-order microscopic cell function $N_1(\bm{y})$, higher-order microscopic cell function ${N_{11}}(\bm{y})$ and macroscopic homogenized solution ${u_0}(\bm{x})$ calculated by HOMS-DRM and DNS-FEM when the thermal conductivities $a_{11} = a_{22} = 0.1$ of inclusion material in Figs.\hspace{1mm}\ref{fig:4}, \ref{fig:5} and \ref{fig:6}.

\begin{figure}
	\centering
	\includegraphics[width=130mm]{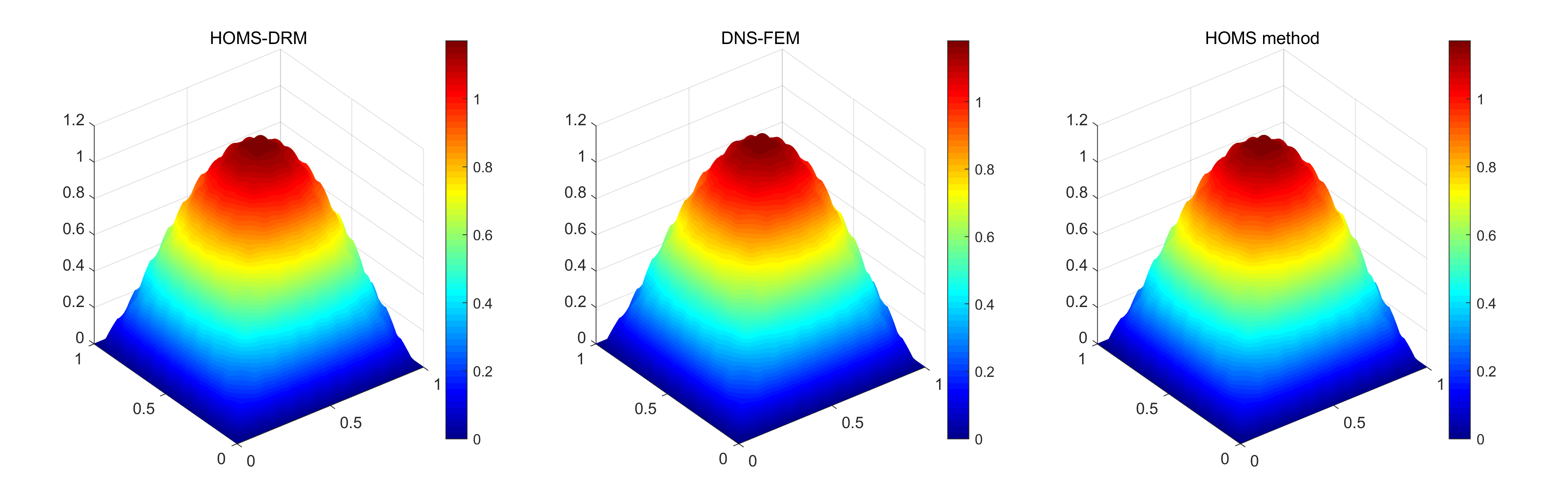}
	\caption{The computational results of temperature field when $a_{11} = a_{22} = 0.1$ of inclusion material.}\label{fig:7}
\end{figure}

\begin{figure}
	\centering
	\includegraphics[width=130mm]{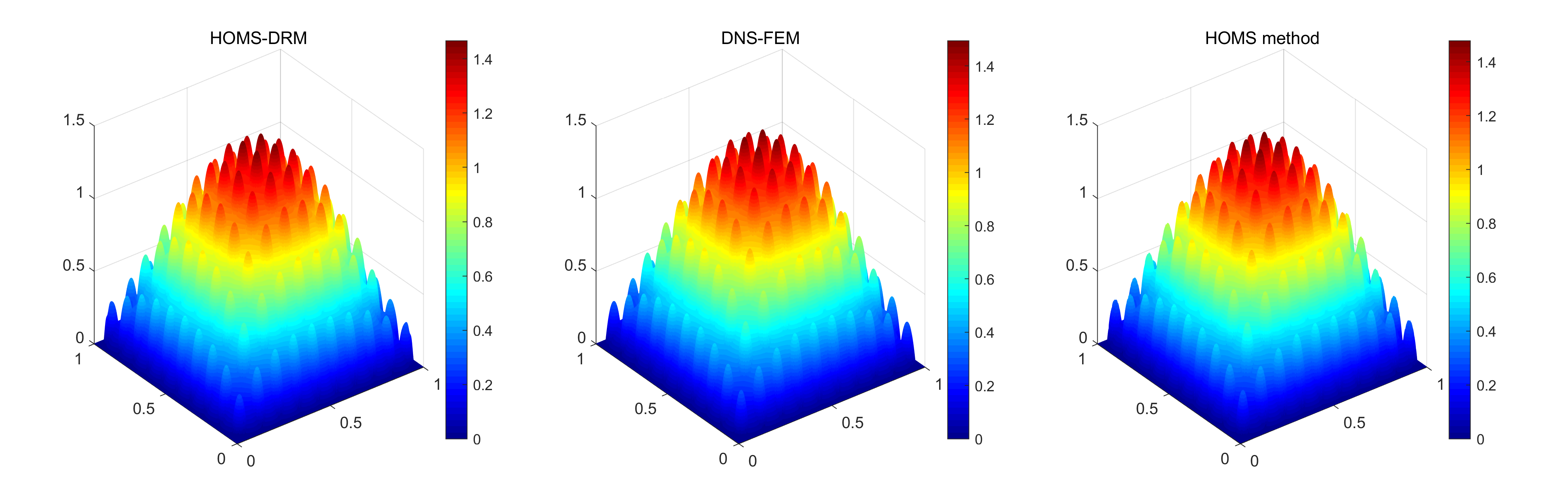}
	\caption{The computational results of temperature field when $a_{11} = a_{22} = 0.01$ of inclusion material.}\label{fig:8}
\end{figure}

\begin{figure}
	\centering
	\includegraphics[width=130mm]{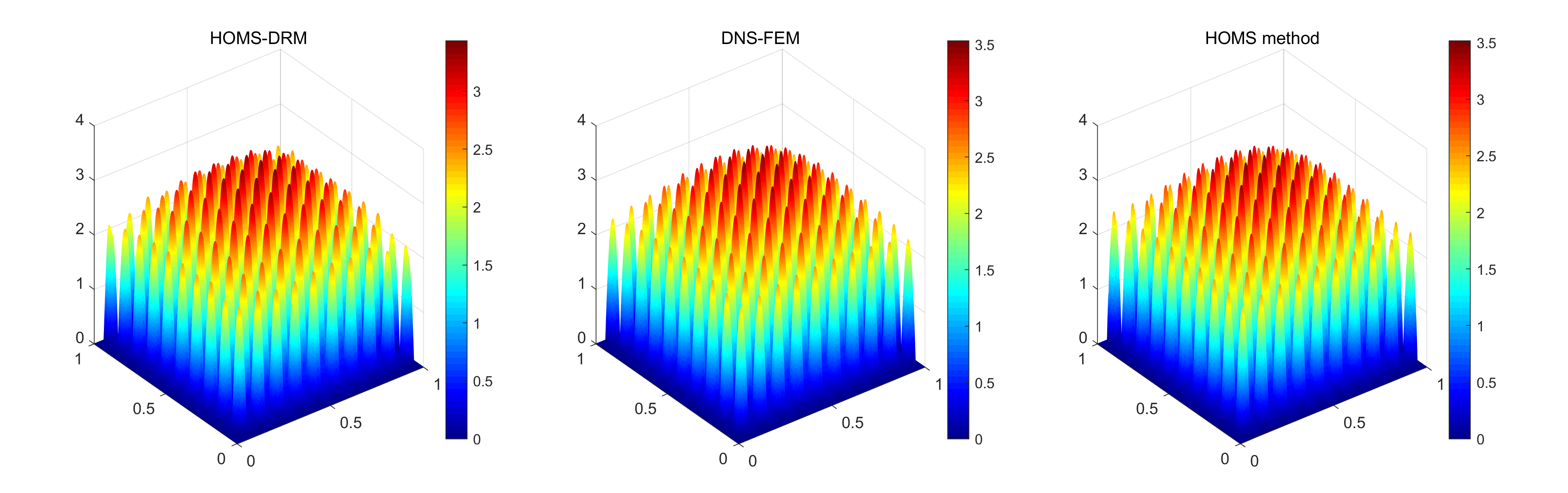}
	\caption{The computational results of temperature field when $a_{11} = a_{22} = 0.001$ of inclusion material.}\label{fig:9}
\end{figure}

Furthermore, Figs.\hspace{1mm}\ref{fig:7}, \ref{fig:8} and \ref{fig:9} display the computational results of HOMS-DRM, DNS-FEM, and HOMS method for simulating steady-state thermal transfer problems of the composites when the thermal conductivities $a_{11} = a_{22}$ of inclusion material are defined as 0.1, 0.01 and 0.001, respectively.

Beside, the specific error results of HOMS-DRM, DNS-FEM and HOMS method are shown in Table \ref{tab:1}.
\begin{table}[]
	
	\centering
	\caption{The error comparisons of thermal transfer simulation.}\label{tab:1}
	\begin{tabular}{cccc}
		\hline
		& $a_{11} = a_{22} =0.1$  & $a_{11} = a_{22} =0.01$ & $a_{11} = a_{22} =0.001$  \\ \hline
		$error_1$ & 0.0104 & 0.0190 & 0.0234 \\
		$error_2$ & 0.0025 & 0.0076 & 0.0218 \\
		$error_3$ & 0.0091 & 0.0129 & 0.0115 \\
        \hline
	\end{tabular}
\end{table}

From the above Figs.\hspace{1mm}\ref{fig:4}, \ref{fig:5} and \ref{fig:6}, firstly we can observe that the proposed HOMS-DRM can accurately solve lower-order and higher-order microscopic cell functions, as well as macroscopic homogenized solution. The relative error ${error_1}$ of lower-order microscopic cell functions ${N_1}(\bm{y})$ and ${N_2}(\bm{y})$ are 0.0107 and 0.0120, respectively. The relative error ${error_1}$ of higher-order microscopic cell functions ${N_{11}}(\bm{y})$, ${N_{12}}(\bm{y})$, ${N_{22}}(\bm{y})$ and ${N_{21}}(\bm{y})$ are 0.0133, 0.0290, 0.0129 and 0.0279, respectively. The relative error of macroscopic homogenized coefficient $a_{ij}^0$ is 0.0015, and the relative error ${error_1}$ of macroscopic homogenized solution ${u_0}(\bm{x})$ is 0.0019. It should be noted that accurate computation of microscopic cell functions and macroscopic homogenized solution of composite material is critical to assemble high-accuracy and high-resolution multi-scale solution according to higher-order multi-scale model \ref{eq:4}. Furthermore, it can be observed from Figs.\hspace{1mm}\ref{fig:7}, \ref{fig:8} and \ref{fig:9} that the proposed HOMS-DRM in this study can accurately simulate multi-scale steady-state thermal transfer problems and capture the microscopic oscillatory information of composite materials, especially for high-contrast composite materials when $a_{11} = a_{22} =0.001$ of inclusion material. Moreover, according to the error results in Table 1, we can conclude that the numerical solutions of HOMS-DRM are very close to those of DNS-FEM, and HOMS method. The presented HOMS-DRM is a potential alternative to traditional numerical-driven HOMS method and also can overcome limitations of prohibitive computation and Frequency Principle when direct deep learning simulation for multi-scale thermal transfer problems.
\subsection{Example 2: Thermal transfer simulation of composite materials with square inclusions}
\label{subsec:2}
\begin{figure}
	\centering
	\includegraphics[width=110mm]{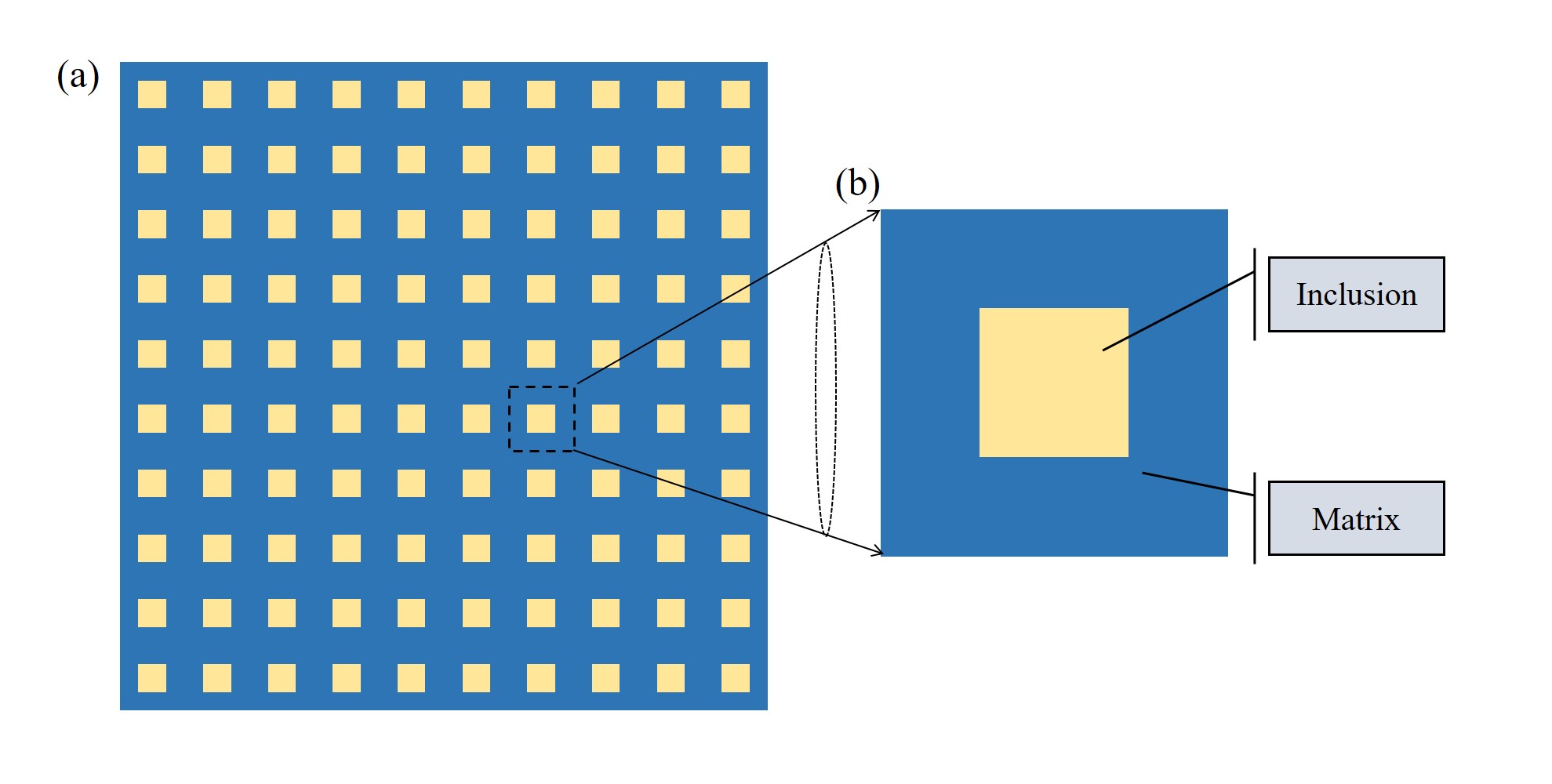}
	\caption{The diagram of periodic composite material. (a) The macroscopic structure $\Omega$; (b) The microscopic unit cell $Y$.}\label{fig:10}
\end{figure}

This example investigates the thermal transfer simulation of composite materials with square inclusions, as detailed in Figs.\hspace{1mm}\ref{fig:10}. The square inclusions of the investigated composite materials have four singular points at their corners, which arise new computational difficulty. At micro-scale, we set $\xi  = 1/10$. For microscopic unit cell $Y$, the side length of square inclusion is 0.5 with central coordinate (0.5,0.5). Moreover, other experimental setup including internal thermal source and material parameters of this example is the same as those of Example 1.

In this example, we utilize a deep residual network with network structure 2-10-[10-10]$\times$4-10-1 to solve the microscopic cell problems and macroscopic homogenized problem. To implement deep learning of HOMS-DRM, 10,000 discrete points are randomly selected based on a uniform distribution within microscopic unit cell $Y$. After 30,000 training epochs with Adam optimizer (initial learning rate 0.005 and dynamic decrease to 0.95 times of previous learning rate after every 500 epochs), we obtain learned solutions for lower-order microscopic cell functions $N_{\alpha_1}(\bm{y})$. Next, macroscopic homogenized thermal conductivities are computed based on 40,000 equidistant points within lower-order microscopic cell functions ${N_1}(\bm{y})$. Furthermore, we continue to obtain learned solutions for higher-order microscopic cell functions $N_{\alpha_1\alpha_2}(\bm{y})$ and macroscopic homogenized solution ${u_0}(\bm{x})$ with 10,000 random points and after 30,000 training epochs.

\begin{figure}
	\centering
	\includegraphics[width=130mm]{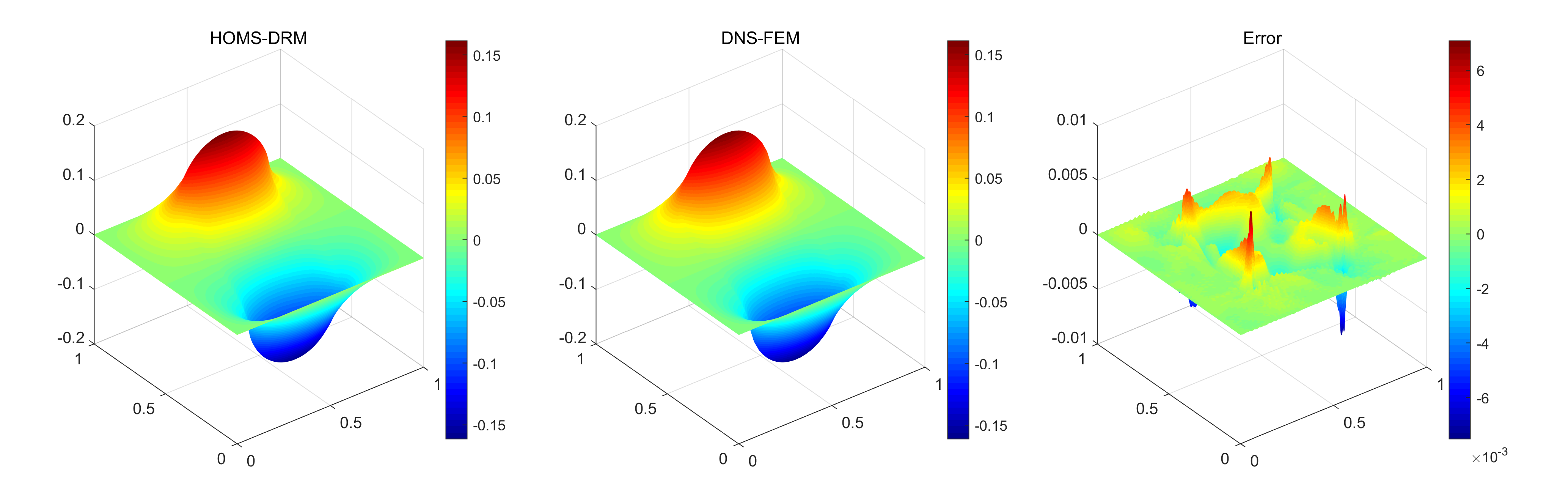}
	\caption{The computational results of lower-order microscopic cell function $N_1(\bm{y})$.}\label{fig:11}
\end{figure}

\begin{figure}
	\centering
	\includegraphics[width=130mm]{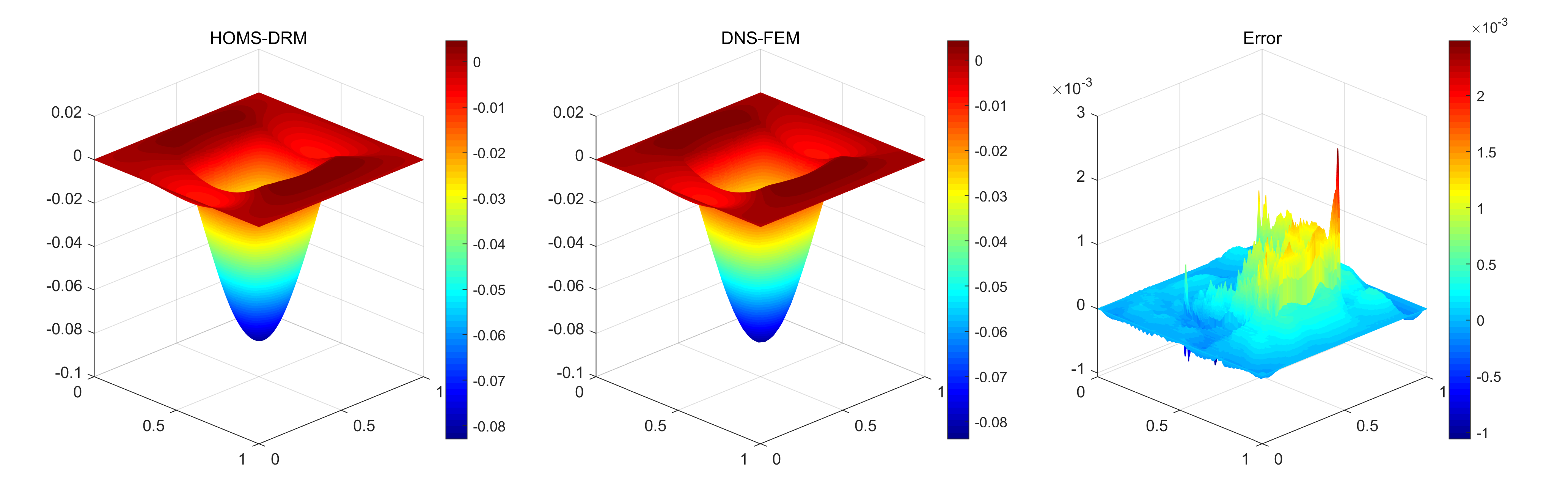}
	\caption{The computational results of higher-order microscopic cell function $N_{11}(\bm{y})$.}\label{fig:12}
\end{figure}

\begin{figure}
	\centering
	\includegraphics[width=130mm]{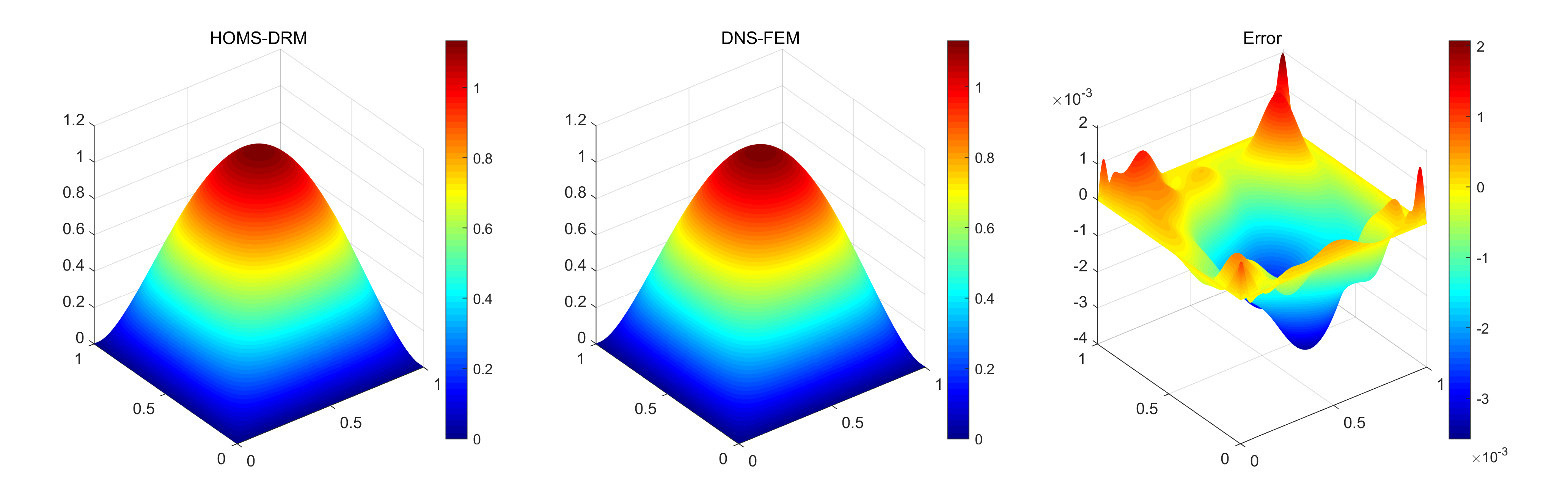}
	\caption{The computational results of macroscopic homogenized solution $u_0(\bm{x})$.}\label{fig:13}
\end{figure}

After the training is stabilized, Figs.\hspace{1mm}\ref{fig:11}, \ref{fig:12} and \ref{fig:13} show the numerical solutions and error comparisons of lower-order microscopic cell functions ${N_1}(\bm{y})$, higher-order microscopic cell functions ${N_{11}}(\bm{y})$ and macroscopic homogenized solution ${u_0}(\bm{x})$ calculated by HOMS-DRM and DNS-FEM when setting $a_{11} = a_{22} = 0.1$.

\begin{figure}
	\centering
	\includegraphics[width=130mm]{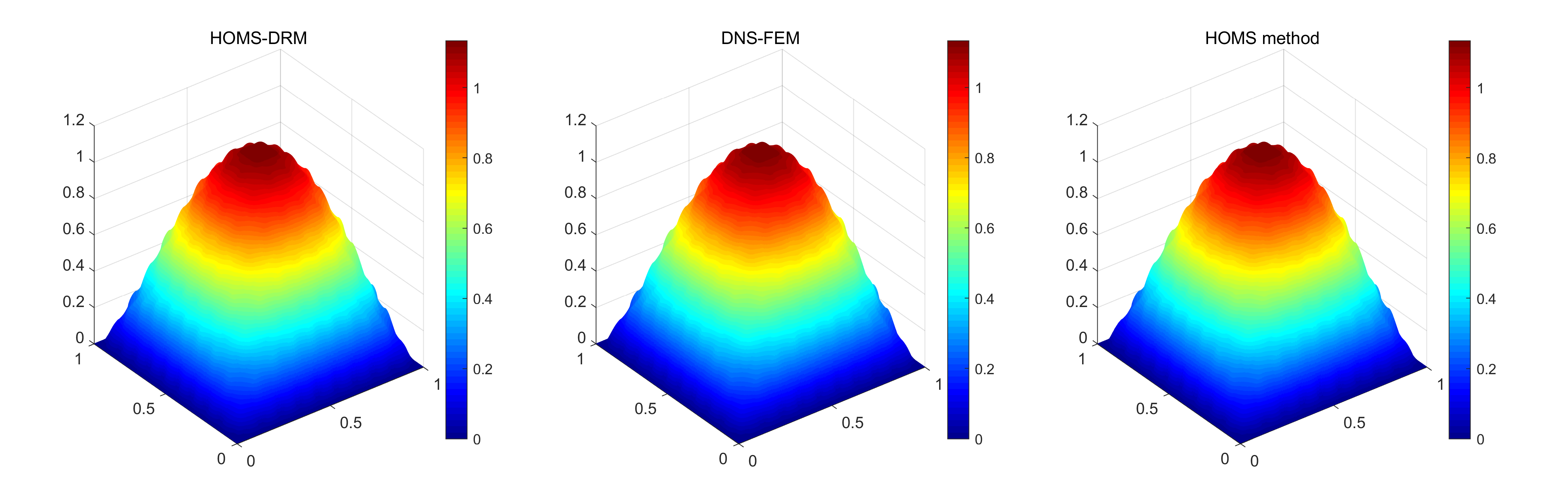}
	\caption{The computational results of temperature field when $a_{11} = a_{22} = 0.1$ of inclusion material.}\label{fig:14}
\end{figure}

\begin{figure}
	\centering
	\includegraphics[width=130mm]{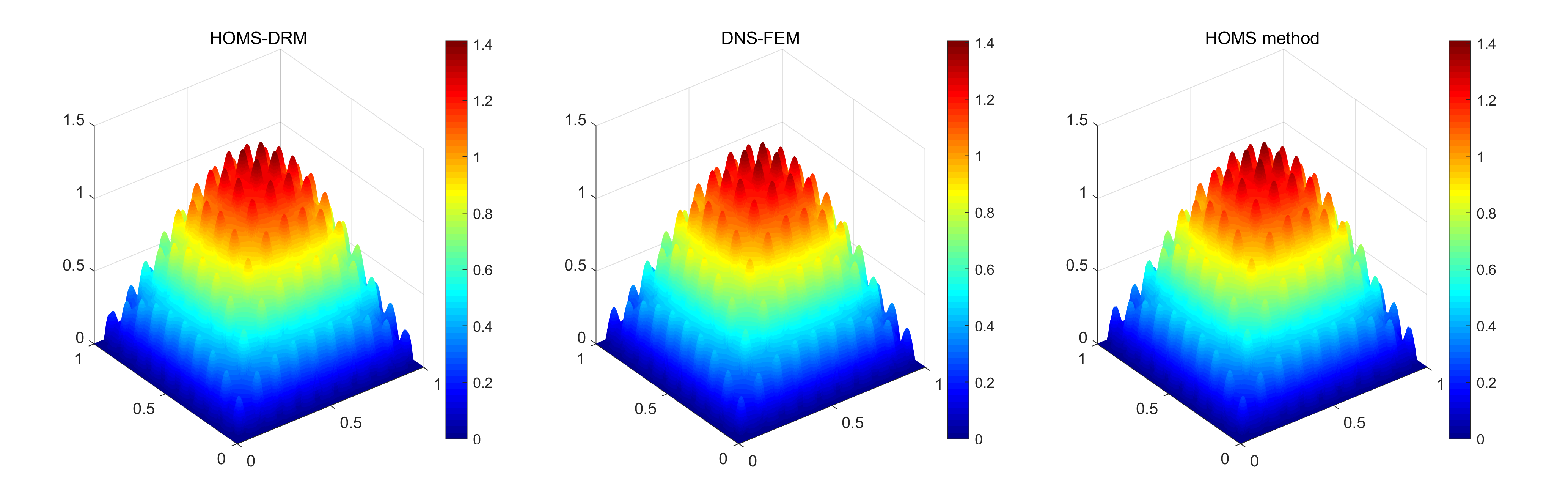}
	\caption{The computational results of temperature field when $a_{11} = a_{22} = 0.01$ of inclusion material.}\label{fig:15}
\end{figure}

\begin{figure}
	\centering
	\includegraphics[width=130mm]{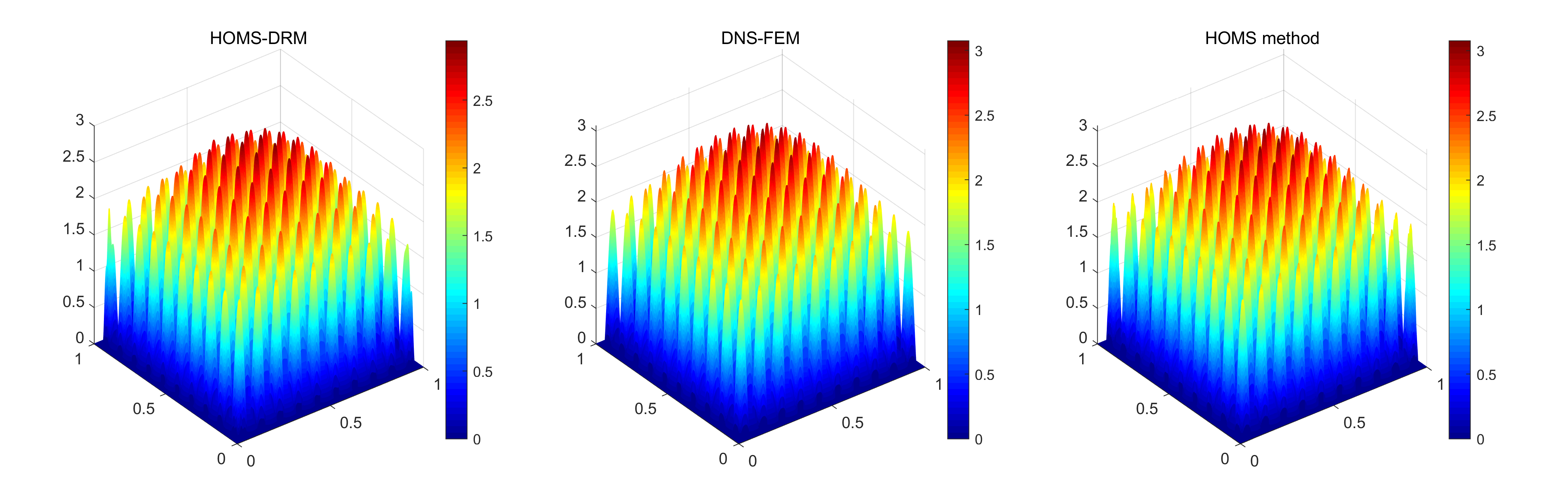}
	\caption{The computational results of temperature field when $a_{11} = a_{22} = 0.001$ of inclusion material.}\label{fig:16}
\end{figure}

The simulative results of steady-state thermal transfer problems of the investigated composites via HOMS-DRM, DNS-FEM and HOMS method are illustrated in Figs.\hspace{1mm}\ref{fig:14}, \ref{fig:15} and \ref{fig:16} for thermal conductivities $a_{11} = a_{22}=0.1$, $0.01$ and $0.001$ of inclusion material, respectively.

To display the computational performance of HOMS-DRM, the specific error comparisons of HOMS-DRM, DNS-FEM and HOMS method are further detailed in Table \ref{tab:2}.
\begin{table}[]	
	\centering
	\caption{The error comparisons of thermal transfer simulation.}\label{tab:2}
	\begin{tabular}{cccc}
		\hline
		& $a_{11} = a_{22} =0.1$  & $a_{11} = a_{22} =0.01$ & $a_{11} = a_{22} =0.001$ \\ \hline
		$error_1$ & 0.0038 & 0.0073 & 0.0340 \\
		$error_2$ & 0.0013 & 0.0039 & 0.0356 \\
		$error_3$ & 0.0038 & 0.0050 & 0.0062 \\ \hline
	\end{tabular}
\end{table}

From the numerical results in Figs.\hspace{1mm}\ref{fig:11}-\ref{fig:16}, we can clearly find that our HOMS-DRM exhibit high-accuracy computing performance for solving microscopic cell functions, macroscopic homogenized solution and assembling higher-order multi-scale solutions for thermal transfer problems of composite materials with square inclusions. The ${error_1}$ of lower-order microscopic cell functions ${N_1}(\bm{y})$ and ${N_2}(\bm{y})$ are 0.0143 and 0.0166, respectively. The ${error_1}$ of higher-order microscopic cell functions ${N_{11}}(\bm{y})$, ${N_{12}}(\bm{y})$, ${N_{22}}(\bm{y})$ and ${N_{21}}(\bm{y})$ are 0.0146, 0.0277, 0.0485 and 0.0227 separately. The relative error of macroscopic homogenized coefficient $a_{ij}^0$ is 0.0008 and the ${error_1}$ of macroscopic homogenized solution ${u_0}(\bm{x})$ is 0.0022. Furthermore, the error comparisons in Table \ref{tab:2} indicate the HOMS-DRM is robust and accurate for simulating multi-scale thermal transfer behaviors of composite materials with geometric singular points, especially in high-contrast case.
\subsection{Example 3: Thermal transfer simulation of composite materials with different volume fraction}
\begin{figure}
	\centering
	\includegraphics[width=110mm]{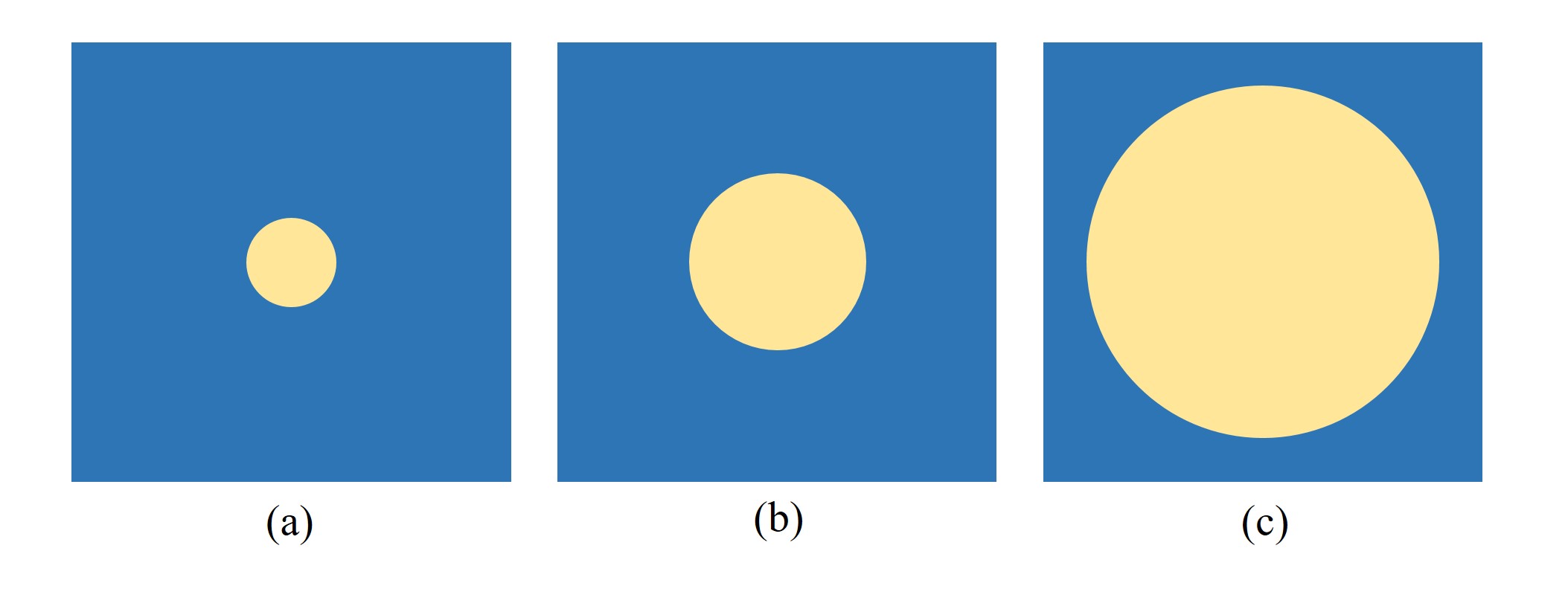}
	\caption{The sketch map of microscopic unit cells with different volume fraction. (a) $r=0.1$; (b) $r=0.2$; (c) $r=0.4$.}\label{fig:17}
\end{figure}

To further validate the computational performance of HOMS-DRM, this example studies thermal transfer simulation of composite materials with different volume fraction, as exhibited in Figs.\hspace{1mm}\ref{fig:17}. The simulated composites all have circular inclusions, whose microscopic unit cells have central coordinate (0.5,0.5) and different radii $r=$0.1, 0.2 and 0.4 of circular inclusions.
\label{subsec:3}

For this example, we still define the internal thermal source $f(\bm{x}) = 10$ and conduct three groups of experiments with $a_{11} = a_{22}=0.1$, $0.01$ and $0.001$ of inclusion material for composite materials possessing fixed volume fraction. Next, we utilize a deep residual network with network structure 2-10-[10-10]$\times$4-10-1 to solve solve the microscopic cell problems and macroscopic homogenized problem. To learn microscopic cell functions and macroscopic homogenized solution, 10,000 computational nodes are randomly selected on the basis of a uniform distribution within microscopic unit cell $Y$ and macroscopic structure $\Omega$. Ending up with 30,000 training epochs with Adam optimizer (initial learning rate 0.005 and dynamic decrease to 0.95 times of previous learning rate after every 500 epochs), the deep learning solutions for microscopic cell problems and macroscopic homogenized problem are gained. Also, macroscopic homogenized thermal parameters are evaluated on a uniform grid with 40,000 nodes within lower-order microscopic cell functions.

After accomplishing deep learning, the evaluation errors of HOMS-DRM, DNS-FEM and HOMS method for composite materials with different volume fraction are given in Tables \ref{tabl:3}, \ref{tabl:4} and \ref{tabl:5}.
\begin{table}[]
	\centering
	\caption{The error comparisons when $a_{11} = a_{22} =0.1$ of inclusion material.}\label{tabl:3}
	\begin{tabular}{ccccc}
		\hline
		& $r=0.1$  & $r=0.2$  & $r=0.3$ & $r=0.4$ \\ \hline
		$error_1$ & 0.0019 & 0.0060 & 0.0104 & 0.0411 \\
		$error_2$ & 0.0026 & 0.0057 & 0.0025 & 0.0017 \\
		$error_3$ & 0.0010 & 0.0024 & 0.0091 & 0.0421 \\ \hline
	\end{tabular}
\end{table}
\begin{table}[]
	\centering
	\caption{The error comparisons when $a_{11} = a_{22} =0.01$ of inclusion material.}\label{tabl:4}
	\begin{tabular}{ccccc}
		\hline
		& $r=0.1$  & $r=0.2$  & $r=0.3$ & $r=0.4$ \\ \hline
		$error_1$ & 0.0026 & 0.0046 & 0.0190 & 0.0622 \\
		$error_2$ & 0.0016 & 0.0037 & 0.0076 & 0.0076 \\
		$error_3$ & 0.0012 & 0.0028 & 0.0129 & 0.0646 \\ \hline
	\end{tabular}
\end{table}
\begin{table}[]
	\centering
	\caption{The error comparisons when $a_{11} = a_{22} =0.001$ of inclusion material.}\label{tabl:5}
	\begin{tabular}{ccccc}
		\hline
		& $r=0.1$  & $r=0.2$  & $r=0.3$ & $r=0.4$ \\ \hline
		$error_1$ & 0.0176 & 0.0323 & 0.0234 & 0.0638 \\
		$error_2$ & 0.0175 & 0.0317 & 0.0218 & 0.0447 \\
		$error_3$ & 0.0011 & 0.0048 & 0.0115 & 0.0385 \\ \hline
	\end{tabular}
\end{table}

\begin{figure}
	\centering
	\includegraphics[width=130mm]{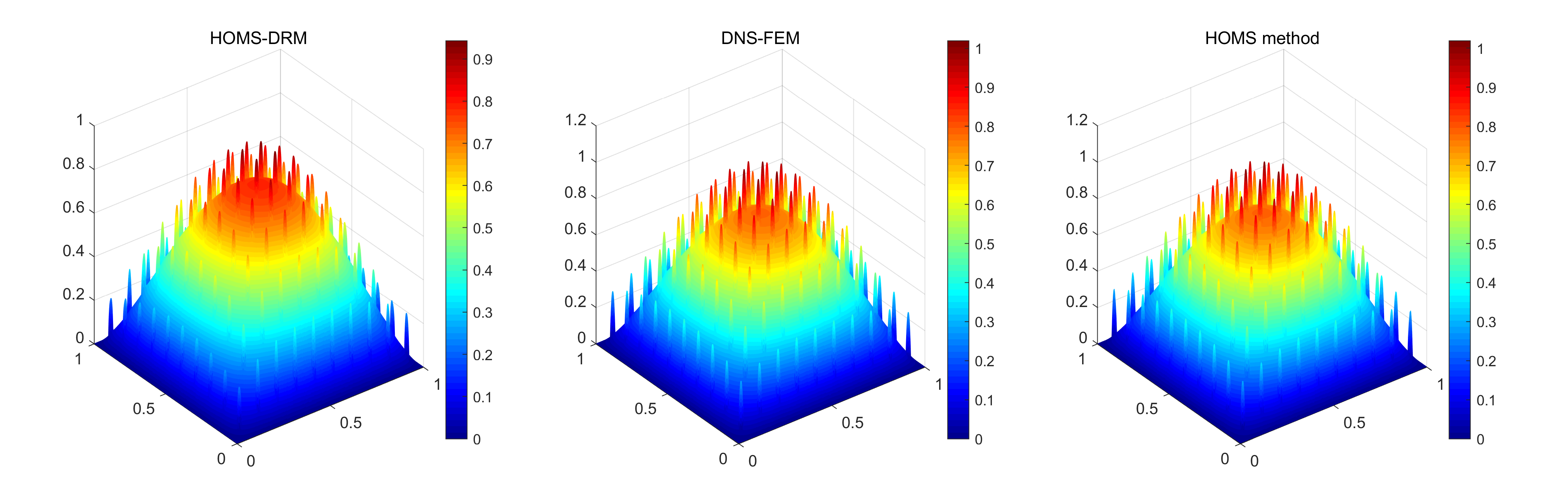}
	\caption{The computational results of temperature field when $r=0.1$.}\label{fig:18}
\end{figure}

\begin{figure}
	\centering
	\includegraphics[width=130mm]{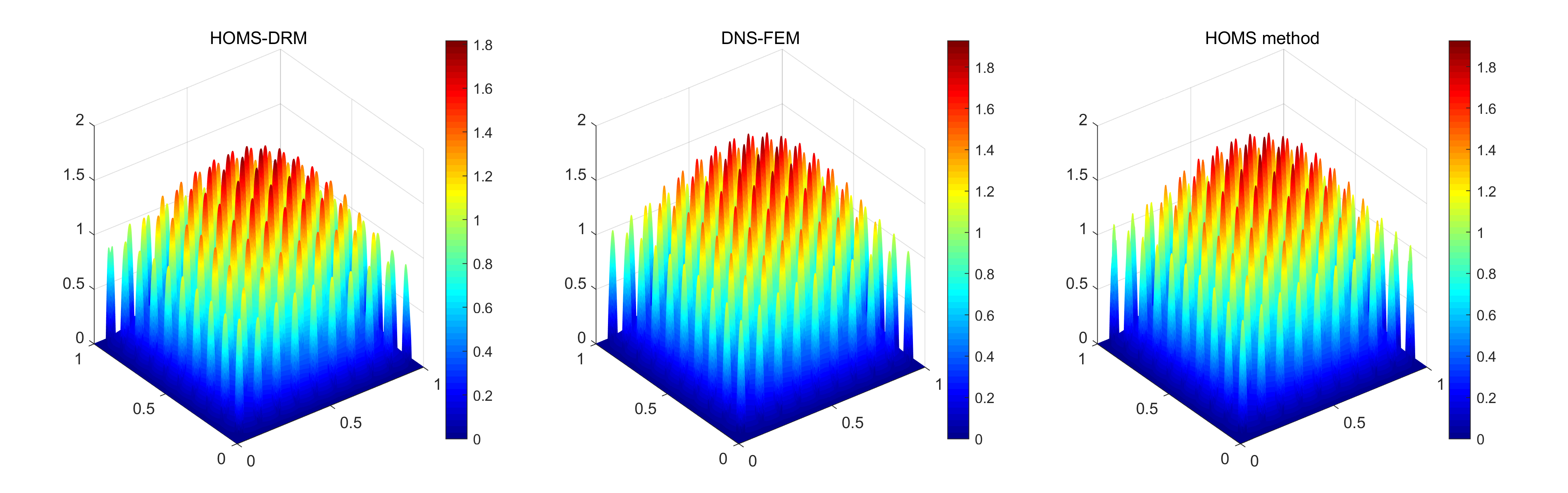}
	\caption{The computational results of temperature field when $r=0.2$.}\label{fig:19}
\end{figure}

\begin{figure}
	\centering
	\includegraphics[width=130mm]{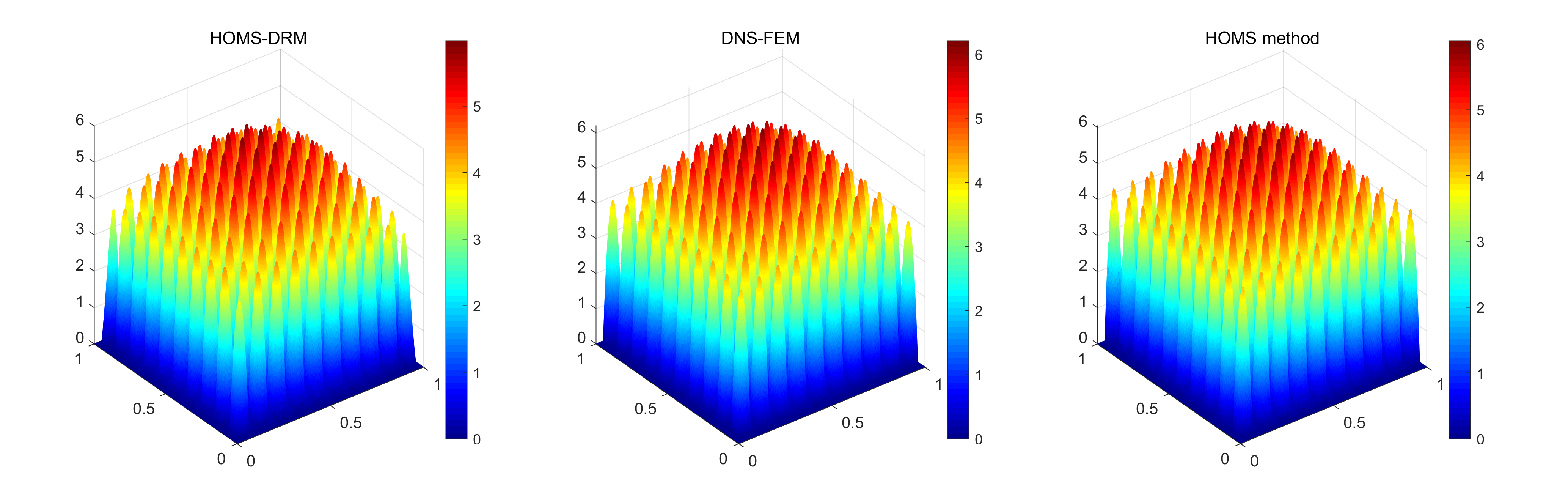}
	\caption{The computational results of temperature field when $r=0.4$.}\label{fig:20}
\end{figure}

Furthermore, some representative cases are exhibited in Figs.\hspace{1mm}\ref{fig:18}, \ref{fig:19} and \ref{fig:20} respectively, which are high-contrast composite materials with the thermal conductivities $a_{11} = a_{22} = 0.001$ of inclusion material.

In a summary, observing the error comparisons in Tables \ref{tabl:3}, \ref{tabl:4} and \ref{tabl:5} can find that the numerical errors of HOMS-DRM simulation will slightly increase when the volume fraction of inclusion material improves. In addition, it can be clearly seen from the above Figs.\hspace{1mm}\ref{fig:18}, \ref{fig:19} and \ref{fig:20} that the oscillatory phenomenon of composite materials is more obvious. Nevertheless, the presented HOMS-DRM can effectively simulate the temperature field of composite materials with different volume fraction, which is of great practical values for engineering simulation.
\subsection{Example 4: Thermal transfer simulation of composite materials with multiple inclusions}
\label{subsec:4}
\begin{figure}
	\centering
	\includegraphics[width=110mm]{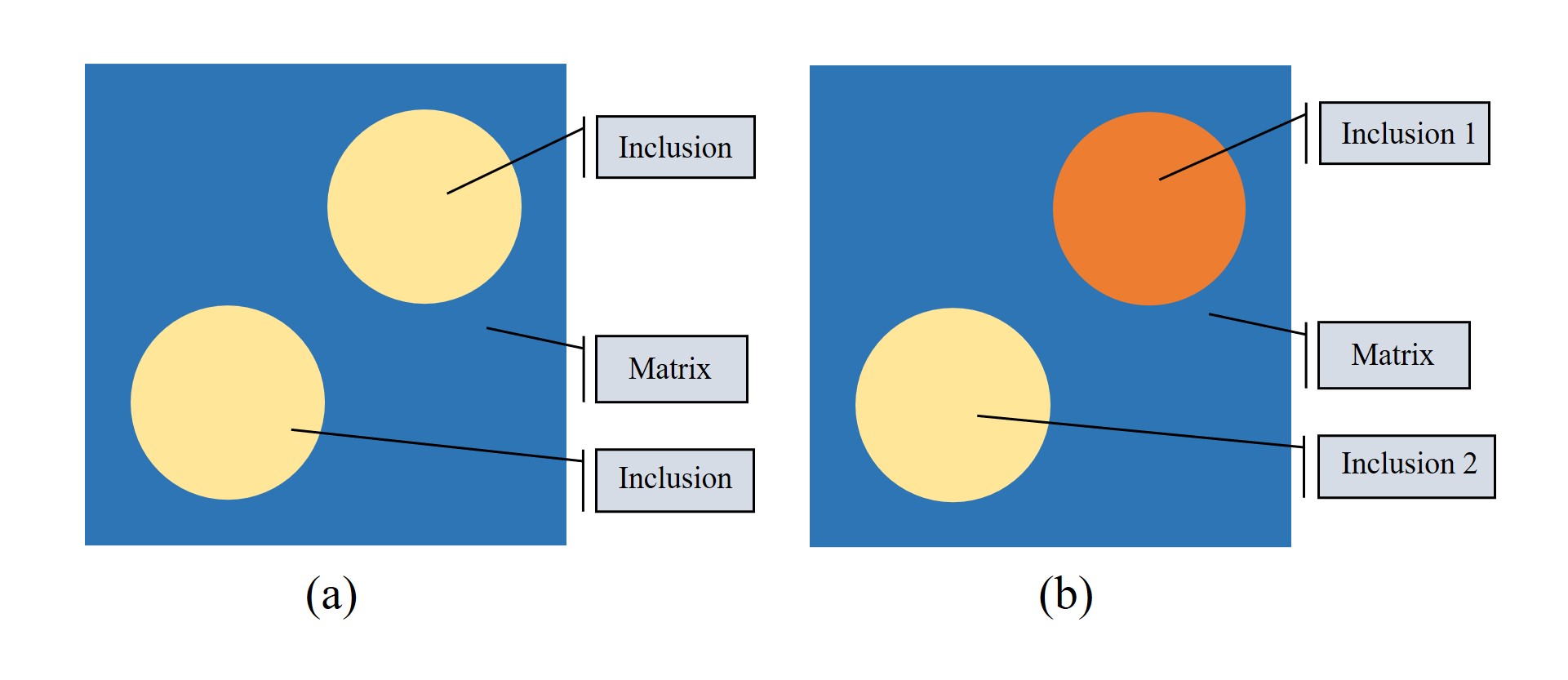}
	\caption{The sketch map of microscopic unit cells with multiple inclusions. (a) Multiple inclusions with same constituent; (b) Multiple inclusions with distinct constituent.}\label{fig:21}
\end{figure}

Consider composite materials which consists of multiple inclusions embedded in matrix, as shown in Fig.\hspace{1mm}\ref{fig:21}. In terms of the geometry of the composites, the central coordinates of two circular inclusions are set as (0.3,0.3) and (0.7,0.7) respectively, and the radius of them is 0.2. In terms of the material property of the composites, for the composite in Fig.\hspace{1mm}\ref{fig:21}(a), its thermal transfer coefficient of matrix material is set as 1, and the thermal transfer coefficients of two inclusions all are set as 0.01. For the composite in Fig.\hspace{1mm}(b), we define the thermal transfer coefficients in material component, inclusion component with circle center (0.3,0.3) and  inclusion component with circle center (0.7,0.7) as 1, 0.01 and 0.02, respectively.

Following the same setup of internal thermal source as Examples 1-3, we employ deep residual networks to solve microscopic cell problems and macroscopic homogenized problem of the composite in Fig.\hspace{1mm}\ref{fig:21}(a). For lower-order microscopic cell functions ${N_1}(\bm{y})$ and ${N_2}(\bm{y})$, and macroscopic homogenized solution ${u_0}(\bm{x})$, deep residual networks with network structure 2-10-[10-10]$\times$4-10-1 are trained with 10,000 randomly discrete points in computational domain and stabilized after 30,000 training epochs with Adam optimizer (initial learning rate 0.005 and dynamic decrease to 0.95 times of previous learning rate per 500 epochs). The macroscopic homogenized material coefficients are calculated based on 40,000 equidistant points inside lower-order microscopic cell functions. For higher-order microscopic cell functions ${N_{11}}(\bm{y})$, ${N_{12}}(\bm{y})$, ${N_{21}}(\bm{y})$ and ${N_{22}}(\bm{y})$, deep residual networks with network structure 2-10-[10-10]$\times$6-10-1 are trained with 15,000 randomly discrete points in microscopic unit cell. After 50,000 training epochs with Adam optimizer (initial learning rate 0.0005 and dynamic decrease to 0.95 times of previous learning rate per 1,000 epochs), the learned solutions for higher-order microscopic cell functions are obtained. Then, we employ HOMS-DRM to obtain high-accuracy multi-scale solutions for steady-state thermal transfer problems of composite materials with multiple inclusions.

\begin{figure}
	\centering
	\includegraphics[width=130mm]{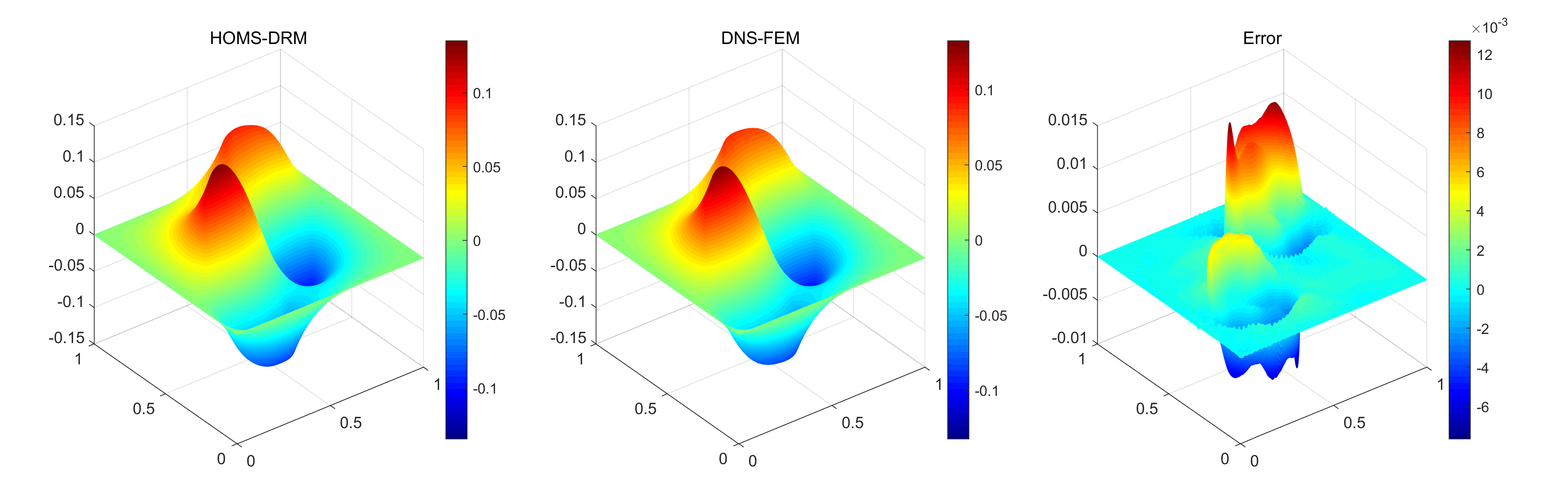}
	\caption{The computational results of $N_1(\bm{y})$.}\label{fig:22}
\end{figure}

\begin{figure}
	\centering
	\includegraphics[width=130mm]{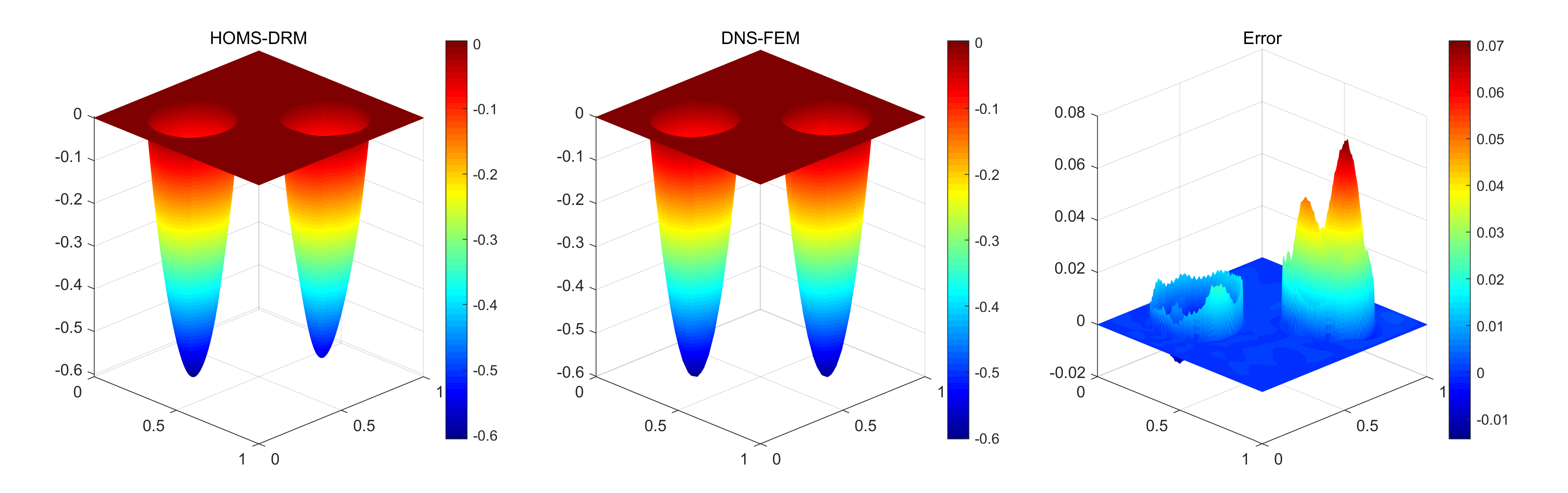}
	\caption{The computational results of $N_{11}(\bm{y})$.}\label{fig:23}
\end{figure}

\begin{figure}
	\centering
	\includegraphics[width=130mm]{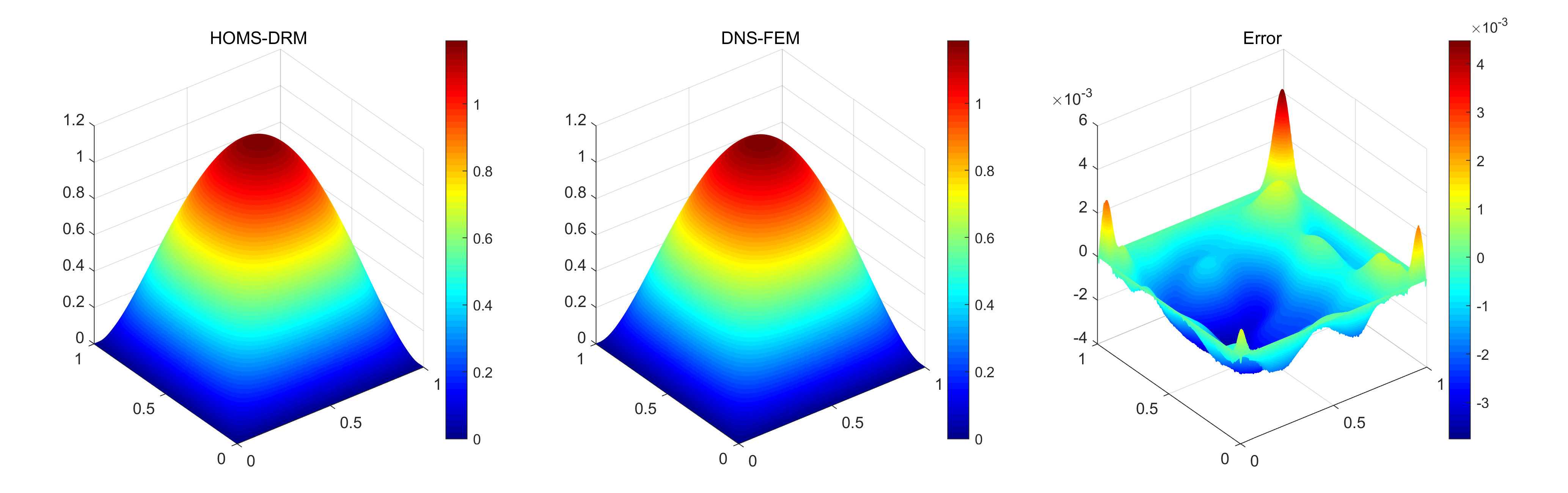}
	\caption{The computational results of $u_0(\bm{x})$.}\label{fig:24}
\end{figure}

\begin{figure}
	\centering
	\includegraphics[width=130mm]{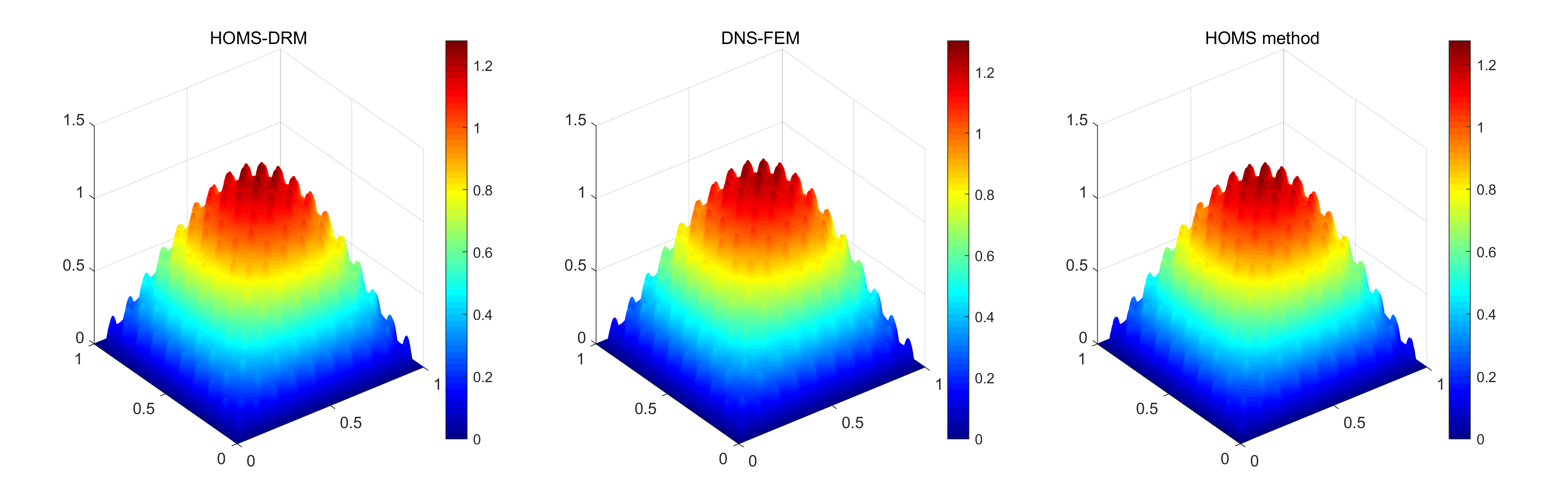}
	\caption{The computational results of temperature field.}\label{fig:25}
\end{figure}

Afterwards, Figs.\hspace{1mm}\ref{fig:22}, \ref{fig:23} and \ref{fig:24} plot the numerical solutions and error comparisons for ${N_1}(\bm{y})$, ${N_{11}}(\bm{y})$ and ${u_0}(\bm{x})$ computed by HOMS-DRM and DNS-FEM, respectively. The ${error_1}$ of lower-order microscopic cell functions ${N_1}(\bm{y})$ and ${N_2}(\bm{y})$ are 0.0525 and 0.0719, respectively. The ${error_1}$ of higher-order microscopic cell functions ${N_{11}}(\bm{y})$, ${N_{12}}(\bm{y})$, ${N_{22}}(\bm{y})$ and ${N_{21}}(\bm{y})$ are 0.0670, 0.0461, 0.0239 and 0.0530, respectively. The relative error of macroscopic homogenized coefficient ${a_{ij}^0}$ is 0.0010, and the ${error_1}$ of macroscopic homogenized solution ${u_0}(\bm{x})$ is 0.0025. Fig.\hspace{1mm}\ref{fig:25} shows the simulative results of HOMS-DRM, DNS-FEM and HOMS method for simulating steady-state thermal transfer behaviors of composite materials, where $erro{r_1}$=0.0194, $erro{r_2}$=0.0039, and $erro{r_3}$=0.0209.

\begin{figure}
	\centering
	\includegraphics[width=130mm]{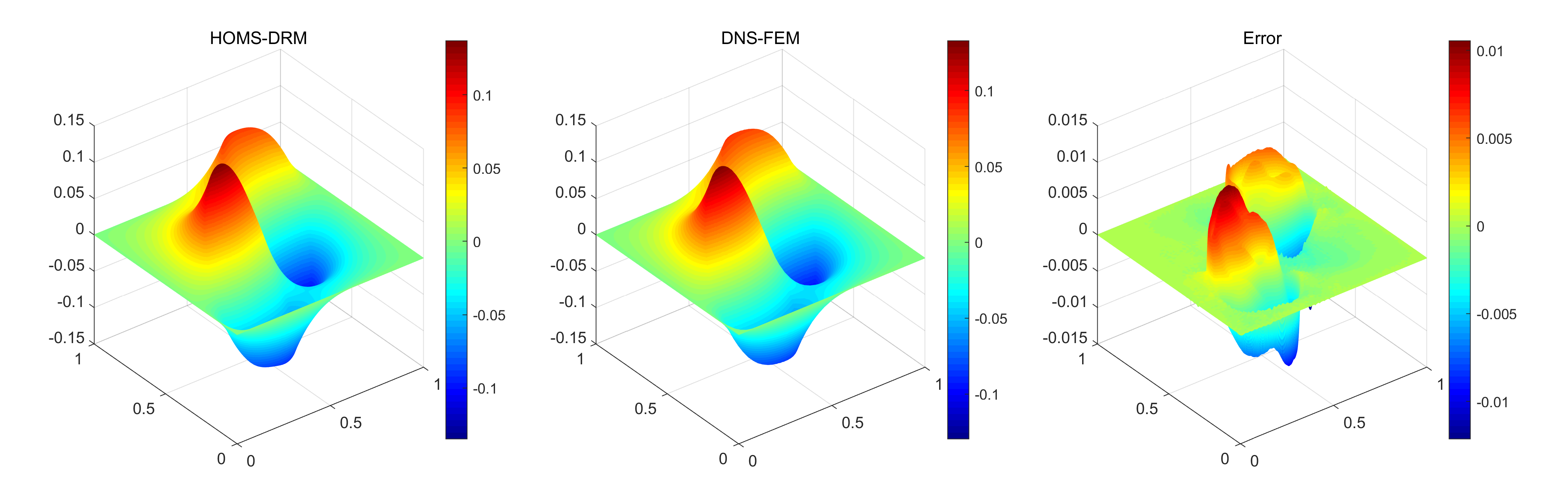}
	\caption{The computational results of $N_1(\bm{y})$.}\label{fig:26}
\end{figure}

\begin{figure}
	\centering
	\includegraphics[width=130mm]{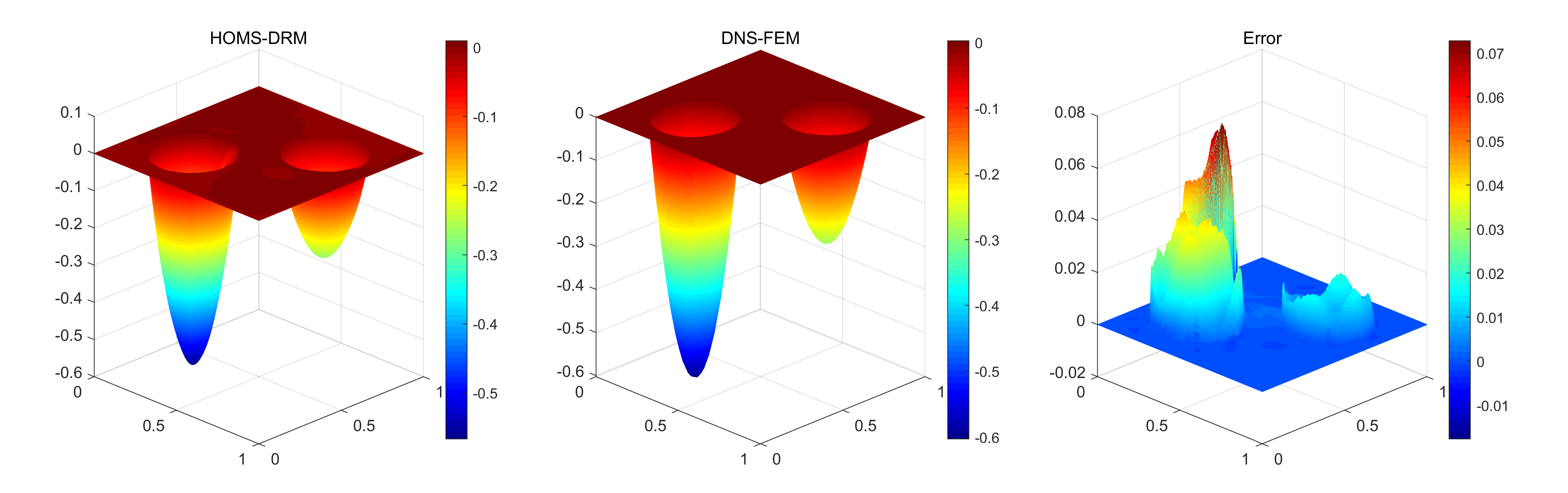}
	\caption{The computational results of $N_{11}(\bm{y})$.}\label{fig:27}
\end{figure}

\begin{figure}
	\centering
	\includegraphics[width=130mm]{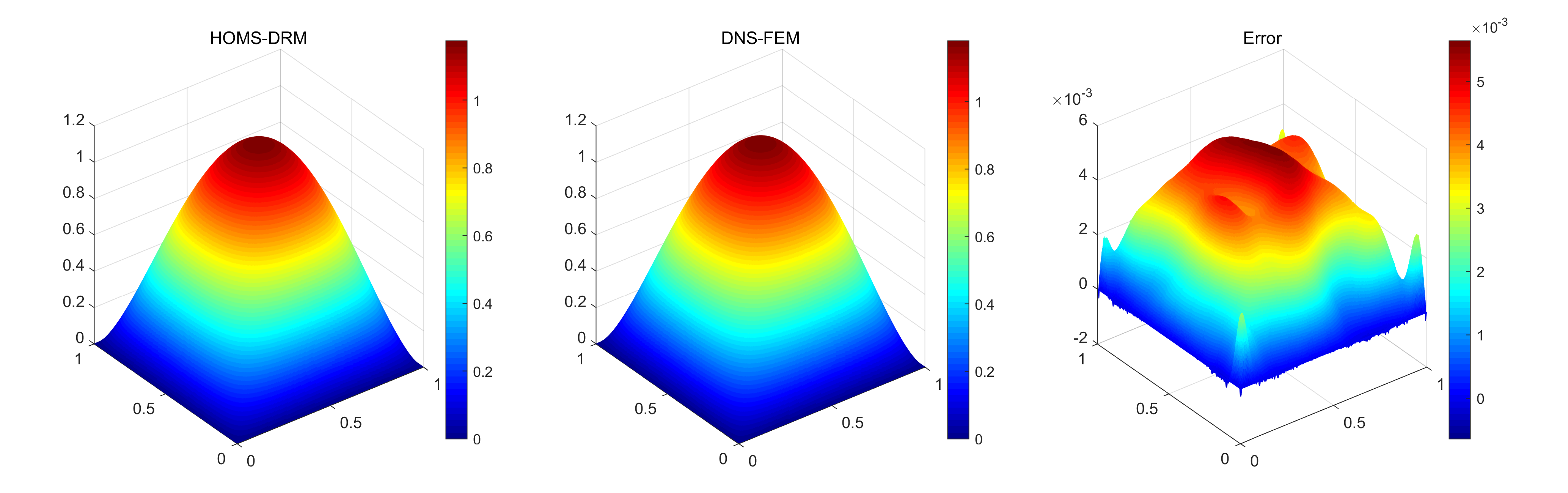}
	\caption{The computational results of $u_0(\bm{x})$.}\label{fig:28}
\end{figure}

\begin{figure}
	\centering
	\includegraphics[width=130mm]{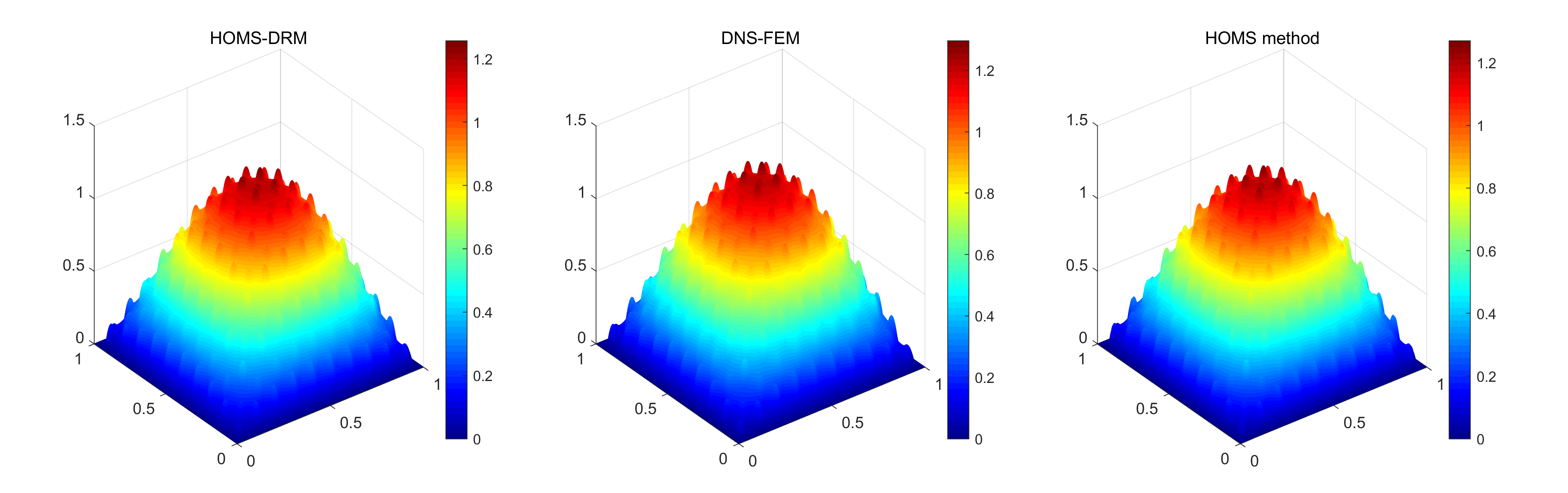}
	\caption{The computational results of temperature field.}\label{fig:29}
\end{figure}

Furthermore, we utilize deep residual networks to solve microscopic cell problems and macroscopic homogenized problem of the composite in Fig.\hspace{1mm}\ref{fig:21}(b). Based on the same setup of learning machine as the composite in Fig.\hspace{1mm}\ref{fig:21}(a), deep learning solutions for lower-order and higher-order microscopic cell functions, macroscopic homogenized solution and higher-order multi-scale solutions are obtained successfully. In the sequel, Figs.\hspace{1mm}\ref{fig:26}, \ref{fig:27} and \ref{fig:28} depict the numerical solutions and error comparisons for ${N_1}(\bm{y})$, ${N_{11}}(\bm{y})$ and ${u_0}(\bm{x})$ computed by HOMS-DRM and DNS-FEM, respectively. The ${error_1}$ of lower-order microscopic cell functions ${N_1}(\bm{y})$ and ${N_2}(\bm{y})$ are 0.0530 and 0.0616 respectively, while ${error_1}$ of higher-order microscopic cell functions ${N_{11}}(\bm{y})$, ${N_{12}}(\bm{y})$, ${N_{22}}(\bm{y})$ and ${N_{21}}(\bm{y})$ are 0.0754, 0.0368, 0.0754 and 0.0320 respectively. In addition, the relative error of macroscopic homogenized coefficient ${a^0}$ is 0.0009, and the ${error_1}$ of macroscopic homogenized solution ${u_0}(\bm{x})$ is 0.0049. Moreover, Fig.\hspace{1mm}\ref{fig:29} shows the computational results of steady-state thermal transfer problems of composite materials via HOMS-DRM, DNS-FEM and HOMS method, where $erro{r_1}$=0.0259, $erro{r_2}$=0.0068, and $erro{r_3}$=0.0205.

In conclusion, the numerical results in Figs.\hspace{1mm}\ref{fig:22}-\ref{fig:29} obviously demonstrate that the proposed HOMS-DRM is effective and reliable to simulate multi-scale heat conduction behaviors of composite materials with complicated multiple inclusions. In addition, highly oscillatory information cased by material heterogeneities are accurately recovered at microscopic level by the proposed HOMS-DRM.
\subsection{Example 5: Thermal transfer simulation of composite materials with different characteristic size of microscopic unit cell}
\label{subsec:5}
It is well-known that the number of microscopic unit cells is $O\big((1/\xi)^2\big)$ in 2D composite materials. In this example, we conduct numerical experiment to evaluate the robustness and stability of our HOMS-DRM by gradually decreasing characteristic size $\xi$. The investigated composite materials have circular inclusions with radius $r=0.3$ in microscopic unit cell. Moreover, the heating source in composite materials is assumed as $f(\bm{x}) = 10$, and high-contrast material  properties of the investigated composites are assumed as $a_{11} = a_{22} =1$ in matrix component and $a_{11} = a_{22} =0.001$ in inclusion component. In this example, we define four different cases as follows. Case I: $\xi=1/5$; Case II: $\xi=1/10$; Case III: $\xi=1/20$; Case IV: $\xi=1/30$.

\begin{figure}
	\centering
	\includegraphics[width=130mm]{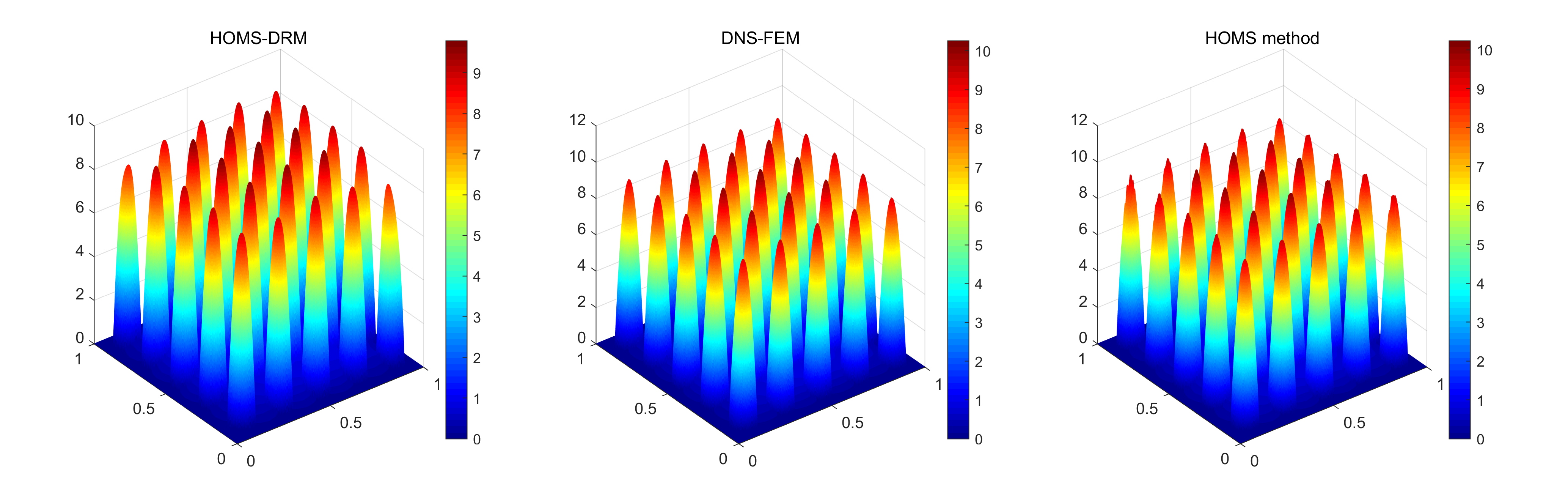}
	\caption{The computational results of temperature field at Case I.}\label{fig:30}
\end{figure}

\begin{figure}
	\centering
	\includegraphics[width=130mm]{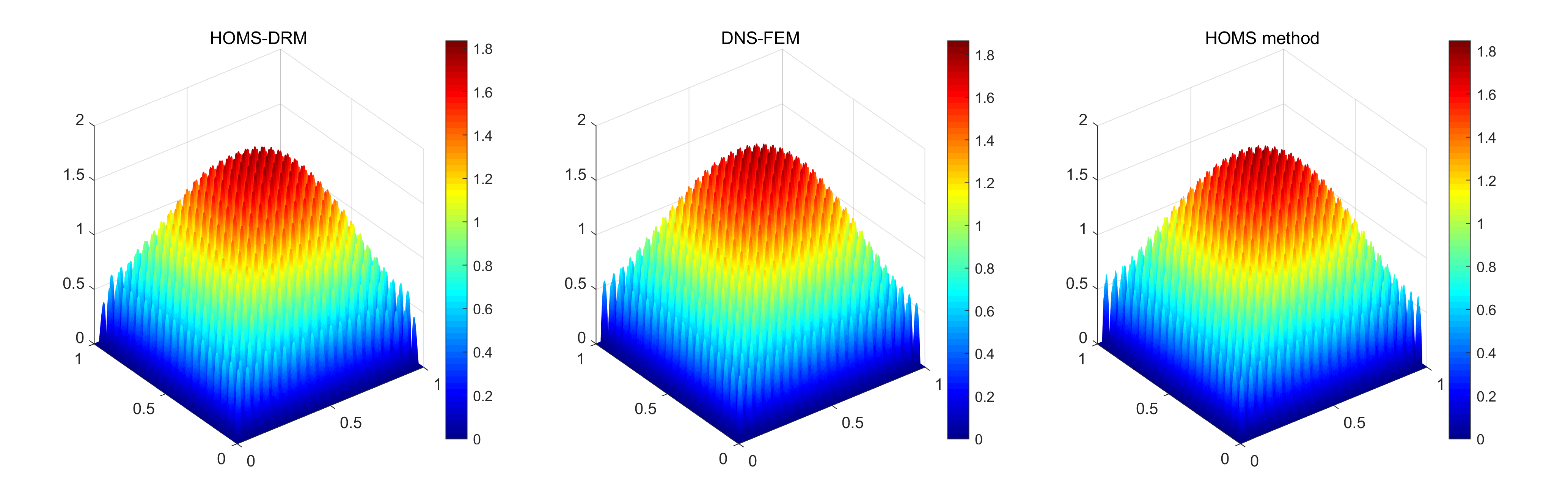}
	\caption{The computational results of temperature field at Case III.}\label{fig:31}
\end{figure}

\begin{figure}
	\centering
	\includegraphics[width=130mm]{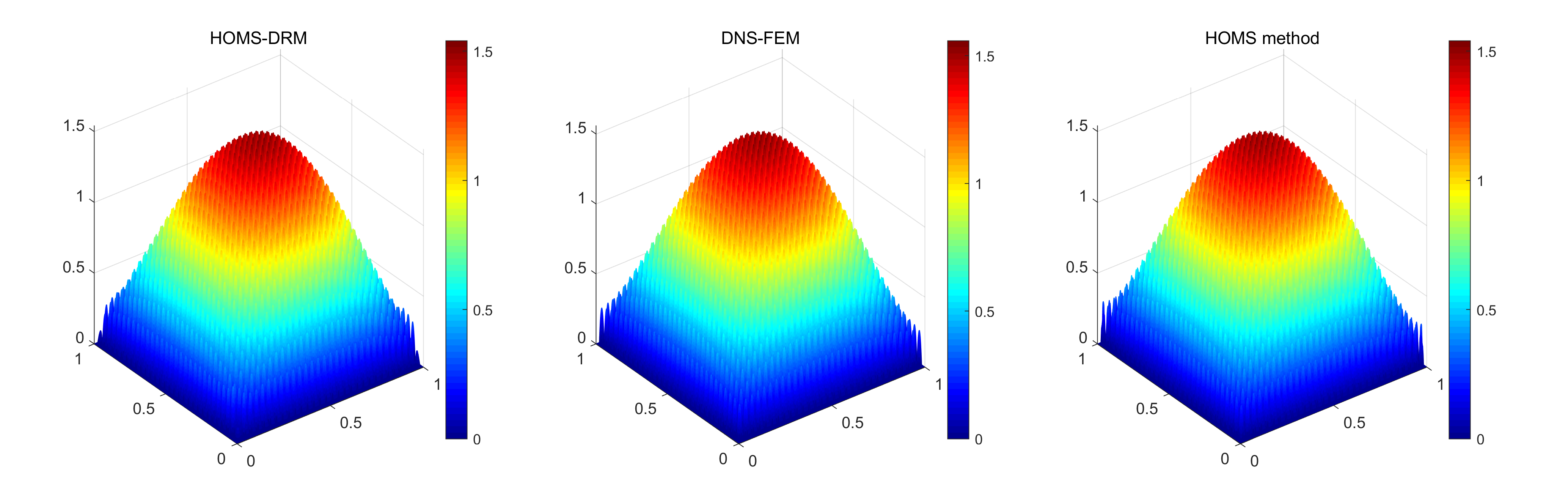}
	\caption{The computational results of temperature field at Case IV.}\label{fig:32}
\end{figure}

Next, deep learning solutions for microscopic cell functions and macroscopic homogenized solution obtained in Example 1 are employed to assemble higher-order multi-scale solutions for thermal transfer simulation of composite materials with the help of automatic differentiation package. Figs.\hspace{1mm}\ref{fig:30}, \ref{fig:31} and \ref{fig:32} depict the simulative results for computing multi-scale temperature field of Cases I, III and IV by HOMS-DRM, DNS-FEM and HOMS method, respectively. The computational results of Case II have been presented in Fig.\hspace{1mm}\ref{fig:9}.

Moreover, the error comparisons of HOMS-DRM, DNS-FEM and HOMS method for thermal transfer simulation of high-contrast material composite materials in four cases are presented in Table \ref{tabl:6}.
\begin{table}[]
	\centering
	\caption{The temperature field calculation error.}\label{tabl:6}
	\begin{tabular}{ccccc}
		\hline
		& $\xi=1/5$  & $\xi=1/10$  & $\xi=1/20$ & $\xi=1/30$ \\ \hline
		$error_1$ & 0.0331 & 0.0234 & 0.0137 & 0.0086 \\
		$error_2$ & 0.0332 & 0.0218 & 0.0139 & 0.0101 \\
		$error_3$ & 0.0077 & 0.0115 & 0.0120 & 0.0128 \\ \hline
	\end{tabular}
\end{table}

As shown in Figs.\hspace{1mm}\ref{fig:30}-\ref{fig:32} and Table \ref{tabl:6}, the proposed HOMS-DRM is robust and stable for simulating large-scale composite materials with massive microscopic unit cells. At this time, DNS-FEM will expend a tremendous amount of computational resources, and even no convergence. Besides, classical deep learning simulation can also need to handle limitations of prohibitive computation and Frequency Principle. However, the prominent advantages of mesh-free, data-free, high numerical accuracy, and low computational overhead of HOMS-DRM enable the effective simulation and analysis of multi-scale thermal transfer behaviors of composite materials, especially for capturing microscopic oscillating phenomenon.
\section{Conclusions and outlook}
\label{sec:6}
Deep learning methods have been widely applied in scientific computing, yet there are relatively few neural network-based methods that can effectively solve multi-scale problems of authentic composite materials. In the present work, we propose a novel higher-order multi-scale deep Ritz method (HOMS-DRM), which inherits the computational advantages of higher-order multi-scale method and deep Ritz method. From the theoretical framework of the presented HOMS-DRM, it can be concluded that our novel approach can overcome limitations of prohibitive computation and Frequency Principle of direct deep learning simulation. And also, the presented HOMS-DRM can tackle the non-differentiable problem at the interface layer and accurately capture the microscopic oscillatory information of of authentic composite materials. Besides, effective numerical algorithm is presented to effectively simulate multi-scale thermal transfer problems of authentic composite materials with highly oscillatory and discontinuous coefficients. In this algorithm, improved deep Ritz method with hard-constraint boundary condition is employed to effectively handle Dirichlet boundary. Besides, a strict convergence proof of the proposed HOMS-DRM is derived under certain assumptions. Numerical results on extensive numerical experiments clearly demonstrated that the presented HOMS-DRM can effectively simulate multi-scale thermal transfer behaviors of composite materials and accurately capture highly oscillatory information of authentic composite materials at the micro-scale, which strongly support the theoretical results of this study. The advanced advantages of mesh-free, data-free, high numerical accuracy and low computational overhead are of great application values for large-scale engineering problems.

In follow-up work, we will rigorously prove the convergence of the presented HOMS-DRM under general conditions. Besides, broad extensions of this higher-order multi-scale deep learning framework are concerned with 3D simulation for mechanics problems, stochastic multi-scale problems and multi-field coupled problems, etc.
\section*{Acknowledgments}
The authors gratefully acknowledge the support of the the National Natural Science Foundation of China (No.\hspace{1mm}12001414), Young Talent Fund of Association for Science and Technology in Shaanxi, China (No.\hspace{1mm}20220506), Young Talent Fund of Association for Science and Technology in Xi'an, China (No.\hspace{1mm}095920221338), the Fundamental Research Funds for the Central Universities (No.\hspace{1mm}KYFZ23020), the National Natural Science Foundation of China (Nos.\hspace{1mm}11971386 and 51739007).

\appendix
\renewcommand{\appendixname}{Appendix~\Alph{section}}

\bibliographystyle{model1-num-names}
\bibliography{paper}

\begin{thebibliography}{56}
\expandafter\ifx\csname natexlab\endcsname\relax\def\natexlab#1{#1}\fi
\providecommand{\url}[1]{\texttt{#1}}
\providecommand{\href}[2]{#2}
\providecommand{\path}[1]{#1}
\providecommand{\DOIprefix}{doi:}
\providecommand{\ArXivprefix}{arXiv:}
\providecommand{\URLprefix}{URL: }
\providecommand{\Pubmedprefix}{pmid:}
\providecommand{\doi}[1]{\href{http://dx.doi.org/#1}{\path{#1}}}
\providecommand{\Pubmed}[1]{\href{pmid:#1}{\path{#1}}}
\providecommand{\bibinfo}[2]{#2}
\ifx\xfnm\relax \def\xfnm[#1]{\unskip,\space#1}\fi
\bibitem[{Chung et~al.(2016)Chung, Efendiev, and Hou}]{R1}
\bibinfo{author}{E.~Chung}, \bibinfo{author}{Y.~Efendiev},
  \bibinfo{author}{T.~Y. Hou},
\newblock \bibinfo{title}{Adaptive multiscale model reduction with generalized
  multiscale finite element methods},
\newblock \bibinfo{journal}{Journal of Computational Physics}
  \bibinfo{volume}{320} (\bibinfo{year}{2016}) \bibinfo{pages}{69--95}.
\bibitem[{Dong et~al.(2018)Dong, Cui, Nie, and Yang}]{R45}
\bibinfo{author}{H.~Dong}, \bibinfo{author}{J.~Cui}, \bibinfo{author}{Y.~Nie},
  \bibinfo{author}{Z.~Yang},
\newblock \bibinfo{title}{Second-order two-scale computational method for
  damped dynamic thermo-mechanical problems of quasi-periodic composite
  materials},
\newblock \bibinfo{journal}{Journal of Computational and Applied Mathematics}
  \bibinfo{volume}{343} (\bibinfo{year}{2018}) \bibinfo{pages}{575--601}.
\bibitem[{Kalina et~al.(2023)Kalina, Linden, Brummund, and K{\"a}stner}]{R46}
\bibinfo{author}{K.~A. Kalina}, \bibinfo{author}{L.~Linden},
  \bibinfo{author}{J.~Brummund}, \bibinfo{author}{M.~K{\"a}stner},
\newblock \bibinfo{title}{Fe$^{\text{ann}}$: an efficient data-driven
  multiscale approach based on physics-constrained neural networks and
  automated data mining},
\newblock \bibinfo{journal}{Computational Mechanics} \bibinfo{volume}{71}
  (\bibinfo{year}{2023}) \bibinfo{pages}{827--851}.
\bibitem[{Bensoussan et~al.(2011)Bensoussan, Lions, and Papanicolaou}]{R2}
\bibinfo{author}{A.~Bensoussan}, \bibinfo{author}{J.-L. Lions},
  \bibinfo{author}{G.~Papanicolaou}, \bibinfo{title}{Asymptotic analysis for
  periodic structures}, volume \bibinfo{volume}{374},
  \bibinfo{publisher}{American Mathematical Soc.}, \bibinfo{year}{2011}.
\bibitem[{Ole{\"\i}nik et~al.(2009)Ole{\"\i}nik, Shamaev, and Yosifian}]{R7}
\bibinfo{author}{O.~A. Ole{\"\i}nik}, \bibinfo{author}{A.~Shamaev},
  \bibinfo{author}{G.~Yosifian}, \bibinfo{title}{Mathematical problems in
  elasticity and homogenization}, \bibinfo{publisher}{Elsevier},
  \bibinfo{year}{2009}.
\bibitem[{Zhang et~al.(2005)Zhang, Zhang, Guo, and Bi}]{R12}
\bibinfo{author}{H.~Zhang}, \bibinfo{author}{S.~Zhang},
  \bibinfo{author}{X.~Guo}, \bibinfo{author}{J.~Bi},
\newblock \bibinfo{title}{Multiple spatial and temporal scales method for
  numerical simulation of non-classical heat conduction problems: one
  dimensional case},
\newblock \bibinfo{journal}{International Journal of Solids and Structures}
  \bibinfo{volume}{42} (\bibinfo{year}{2005}) \bibinfo{pages}{877--899}.
\bibitem[{Engquist and Souganidis(2008)}]{R3}
\bibinfo{author}{B.~Engquist}, \bibinfo{author}{P.~E. Souganidis},
\newblock \bibinfo{title}{Asymptotic and numerical homogenization},
\newblock \bibinfo{journal}{Acta Numerica} \bibinfo{volume}{17}
  (\bibinfo{year}{2008}) \bibinfo{pages}{147--190}.
\bibitem[{Ming et~al.(2005)Ming, Zhang et~al.}]{R4}
\bibinfo{author}{P.~Ming}, \bibinfo{author}{P.~Zhang}, et~al.,
\newblock \bibinfo{title}{Analysis of the heterogeneous multiscale method for
  elliptic homogenization problems},
\newblock \bibinfo{journal}{Journal of the American Mathematical Society}
  \bibinfo{volume}{18} (\bibinfo{year}{2005}) \bibinfo{pages}{121--156}.
\bibitem[{Hughes et~al.(1998)Hughes, Feij{\'o}o, Mazzei, and Quincy}]{R5}
\bibinfo{author}{T.~J. Hughes}, \bibinfo{author}{G.~R. Feij{\'o}o},
  \bibinfo{author}{L.~Mazzei}, \bibinfo{author}{J.-B. Quincy},
\newblock \bibinfo{title}{The variational multiscale method—a paradigm for
  computational mechanics},
\newblock \bibinfo{journal}{Computer methods in applied mechanics and
  engineering} \bibinfo{volume}{166} (\bibinfo{year}{1998})
  \bibinfo{pages}{3--24}.
\bibitem[{Hou et~al.(1999)Hou, Wu, and Cai}]{R6}
\bibinfo{author}{T.~Hou}, \bibinfo{author}{X.-H. Wu}, \bibinfo{author}{Z.~Cai},
\newblock \bibinfo{title}{Convergence of a multiscale finite element method for
  elliptic problems with rapidly oscillating coefficients},
\newblock \bibinfo{journal}{Mathematics of computation} \bibinfo{volume}{68}
  (\bibinfo{year}{1999}) \bibinfo{pages}{913--943}.
\bibitem[{Efendiev et~al.(2013)Efendiev, Galvis, and Hou}]{R13}
\bibinfo{author}{Y.~Efendiev}, \bibinfo{author}{J.~Galvis},
  \bibinfo{author}{T.~Y. Hou},
\newblock \bibinfo{title}{Generalized multiscale finite element methods
  (gmsfem)},
\newblock \bibinfo{journal}{Journal of computational physics}
  \bibinfo{volume}{251} (\bibinfo{year}{2013}) \bibinfo{pages}{116--135}.
\bibitem[{Dong et~al.(2019)Dong, Zheng, Cui, Nie, Yang, and Yang}]{R8}
\bibinfo{author}{H.~Dong}, \bibinfo{author}{X.~Zheng},
  \bibinfo{author}{J.~Cui}, \bibinfo{author}{Y.~Nie},
  \bibinfo{author}{Z.~Yang}, \bibinfo{author}{Z.~Yang},
\newblock \bibinfo{title}{Multiscale computational method for dynamic
  thermo-mechanical problems of composite structures with diverse periodic
  configurations in different subdomains},
\newblock \bibinfo{journal}{Journal of Scientific Computing}
  \bibinfo{volume}{79} (\bibinfo{year}{2019}) \bibinfo{pages}{1630--1666}.
\bibitem[{Cao et~al.(2001)Cao, Cui, Zhu, and Luo}]{R9}
\bibinfo{author}{L.-q. Cao}, \bibinfo{author}{J.-z. Cui},
  \bibinfo{author}{D.-c. Zhu}, \bibinfo{author}{J.-l. Luo},
\newblock \bibinfo{title}{Multiscale finite element method for subdivided
  periodic elastic structures of composite materials},
\newblock \bibinfo{journal}{Journal of Computational Mathematics}
  (\bibinfo{year}{2001}) \bibinfo{pages}{205--212}.
\bibitem[{Cao et~al.(2002)Cao, Cui, and Zhu}]{R10}
\bibinfo{author}{L.-Q. Cao}, \bibinfo{author}{J.-Z. Cui},
  \bibinfo{author}{D.-C. Zhu},
\newblock \bibinfo{title}{Multiscale asymptotic analysis and numerical
  simulation for the second order helmholtz equations with rapidly oscillating
  coefficients over general convex domains},
\newblock \bibinfo{journal}{SIAM Journal on Numerical Analysis}
  \bibinfo{volume}{40} (\bibinfo{year}{2002}) \bibinfo{pages}{543--577}.
\bibitem[{Cui and Cao(1999)}]{R11}
\bibinfo{author}{J.-z. Cui}, \bibinfo{author}{L.-q. Cao},
\newblock \bibinfo{title}{Two-scale asympototic analysis methods for a class of
  elliptic boundary value problems with small periodic coefficients},
\newblock \bibinfo{journal}{MATHEMATICA NUMERICA SINICA-CHINESE EDITION}
  \bibinfo{volume}{21} (\bibinfo{year}{1999}) \bibinfo{pages}{19--28}.
\bibitem[{Dong et~al.(2022)Dong, Yang, Guan, and Cui}]{R38}
\bibinfo{author}{H.~Dong}, \bibinfo{author}{Z.~Yang},
  \bibinfo{author}{X.~Guan}, \bibinfo{author}{J.~Cui},
\newblock \bibinfo{title}{Stochastic higher-order three-scale strength
  prediction model for composite structures with micromechanical analysis},
\newblock \bibinfo{journal}{Journal of Computational Physics}
  \bibinfo{volume}{465} (\bibinfo{year}{2022}) \bibinfo{pages}{111352}.
\bibitem[{Dong et~al.(2023)Dong, Cui, Nie, Ma, Jin, and Huang}]{R39}
\bibinfo{author}{H.~Dong}, \bibinfo{author}{J.~Cui}, \bibinfo{author}{Y.~Nie},
  \bibinfo{author}{R.~Ma}, \bibinfo{author}{K.~Jin},
  \bibinfo{author}{D.~Huang},
\newblock \bibinfo{title}{Multi-scale computational method for nonlinear
  dynamic thermo-mechanical problems of composite materials with
  temperature-dependent properties},
\newblock \bibinfo{journal}{Communications in Nonlinear Science and Numerical
  Simulation} \bibinfo{volume}{118} (\bibinfo{year}{2023})
  \bibinfo{pages}{107000}.
\bibitem[{Dong et~al.(2019)Dong, Zheng, Cui, Nie, Yang, and Ma}]{R40}
\bibinfo{author}{H.~Dong}, \bibinfo{author}{X.~Zheng},
  \bibinfo{author}{J.~Cui}, \bibinfo{author}{Y.~Nie},
  \bibinfo{author}{Z.~Yang}, \bibinfo{author}{Q.~Ma},
\newblock \bibinfo{title}{Multi-scale computational method for dynamic
  thermo-mechanical performance of heterogeneous shell structures with
  orthogonal periodic configurations},
\newblock \bibinfo{journal}{Computer Methods in Applied Mechanics and
  Engineering} \bibinfo{volume}{354} (\bibinfo{year}{2019})
  \bibinfo{pages}{143--180}.
\bibitem[{Rudy et~al.(2017)Rudy, Brunton, Proctor, and Kutz}]{R14}
\bibinfo{author}{S.~H. Rudy}, \bibinfo{author}{S.~L. Brunton},
  \bibinfo{author}{J.~L. Proctor}, \bibinfo{author}{J.~N. Kutz},
\newblock \bibinfo{title}{Data-driven discovery of partial differential
  equations},
\newblock \bibinfo{journal}{Science advances} \bibinfo{volume}{3}
  (\bibinfo{year}{2017}) \bibinfo{pages}{e1602614}.
\bibitem[{Qin et~al.(2019)Qin, Wu, and Xiu}]{R15}
\bibinfo{author}{T.~Qin}, \bibinfo{author}{K.~Wu}, \bibinfo{author}{D.~Xiu},
\newblock \bibinfo{title}{Data driven governing equations approximation using
  deep neural networks},
\newblock \bibinfo{journal}{Journal of Computational Physics}
  \bibinfo{volume}{395} (\bibinfo{year}{2019}) \bibinfo{pages}{620--635}.
\bibitem[{Lagaris et~al.(1998)Lagaris, Likas, and Fotiadis}]{R47}
\bibinfo{author}{I.~E. Lagaris}, \bibinfo{author}{A.~Likas},
  \bibinfo{author}{D.~I. Fotiadis},
\newblock \bibinfo{title}{Artificial neural networks for solving ordinary and
  partial differential equations},
\newblock \bibinfo{journal}{IEEE transactions on neural networks}
  \bibinfo{volume}{9} (\bibinfo{year}{1998}) \bibinfo{pages}{987--1000}.
\bibitem[{Wang et~al.(2020)Wang, Cheung, Chung, Efendiev, and Wang}]{R41}
\bibinfo{author}{Y.~Wang}, \bibinfo{author}{S.~W. Cheung},
  \bibinfo{author}{E.~T. Chung}, \bibinfo{author}{Y.~Efendiev},
  \bibinfo{author}{M.~Wang},
\newblock \bibinfo{title}{Deep multiscale model learning},
\newblock \bibinfo{journal}{Journal of Computational Physics}
  \bibinfo{volume}{406} (\bibinfo{year}{2020}) \bibinfo{pages}{109071}.
\bibitem[{Kiyani et~al.(2022)Kiyani, Silber, Kooshkbaghi, and Karttunen}]{R16}
\bibinfo{author}{E.~Kiyani}, \bibinfo{author}{S.~Silber},
  \bibinfo{author}{M.~Kooshkbaghi}, \bibinfo{author}{M.~Karttunen},
\newblock \bibinfo{title}{Machine-learning-based data-driven discovery of
  nonlinear phase-field dynamics},
\newblock \bibinfo{journal}{Physical Review E} \bibinfo{volume}{106}
  (\bibinfo{year}{2022}) \bibinfo{pages}{065303}.
\bibitem[{Raissi et~al.(2019)Raissi, Perdikaris, and Karniadakis}]{R17}
\bibinfo{author}{M.~Raissi}, \bibinfo{author}{P.~Perdikaris},
  \bibinfo{author}{G.~E. Karniadakis},
\newblock \bibinfo{title}{Physics-informed neural networks: A deep learning
  framework for solving forward and inverse problems involving nonlinear
  partial differential equations},
\newblock \bibinfo{journal}{Journal of Computational physics}
  \bibinfo{volume}{378} (\bibinfo{year}{2019}) \bibinfo{pages}{686--707}.
\bibitem[{Lu et~al.(2021)Lu, Meng, Mao, and Karniadakis}]{R21}
\bibinfo{author}{L.~Lu}, \bibinfo{author}{X.~Meng}, \bibinfo{author}{Z.~Mao},
  \bibinfo{author}{G.~E. Karniadakis},
\newblock \bibinfo{title}{Deepxde: A deep learning library for solving
  differential equations},
\newblock \bibinfo{journal}{SIAM review} \bibinfo{volume}{63}
  (\bibinfo{year}{2021}) \bibinfo{pages}{208--228}.
\bibitem[{Yu et~al.(2022)Yu, Lu, Meng, and Karniadakis}]{R22}
\bibinfo{author}{J.~Yu}, \bibinfo{author}{L.~Lu}, \bibinfo{author}{X.~Meng},
  \bibinfo{author}{G.~E. Karniadakis},
\newblock \bibinfo{title}{Gradient-enhanced physics-informed neural networks
  for forward and inverse pde problems},
\newblock \bibinfo{journal}{Computer Methods in Applied Mechanics and
  Engineering} \bibinfo{volume}{393} (\bibinfo{year}{2022})
  \bibinfo{pages}{114823}.
\bibitem[{Meng et~al.(2020)Meng, Li, Zhang, and Karniadakis}]{R23}
\bibinfo{author}{X.~Meng}, \bibinfo{author}{Z.~Li}, \bibinfo{author}{D.~Zhang},
  \bibinfo{author}{G.~E. Karniadakis},
\newblock \bibinfo{title}{Ppinn: Parareal physics-informed neural network for
  time-dependent pdes},
\newblock \bibinfo{journal}{Computer Methods in Applied Mechanics and
  Engineering} \bibinfo{volume}{370} (\bibinfo{year}{2020})
  \bibinfo{pages}{113250}.
\bibitem[{Jagtap and Karniadakis(2021)}]{R24}
\bibinfo{author}{A.~D. Jagtap}, \bibinfo{author}{G.~E. Karniadakis},
\newblock \bibinfo{title}{Extended physics-informed neural networks (xpinns): A
  generalized space-time domain decomposition based deep learning framework for
  nonlinear partial differential equations},
\newblock in: \bibinfo{booktitle}{AAAI Spring Symposium: MLPS},
  \bibinfo{year}{2021}, pp. \bibinfo{pages}{2002--2041}.
\bibitem[{Jagtap et~al.(2022)Jagtap, Mao, Adams, and Karniadakis}]{R25}
\bibinfo{author}{A.~D. Jagtap}, \bibinfo{author}{Z.~Mao},
  \bibinfo{author}{N.~Adams}, \bibinfo{author}{G.~E. Karniadakis},
\newblock \bibinfo{title}{Physics-informed neural networks for inverse problems
  in supersonic flows},
\newblock \bibinfo{journal}{Journal of Computational Physics}
  \bibinfo{volume}{466} (\bibinfo{year}{2022}) \bibinfo{pages}{111402}.
\bibitem[{Mao et~al.(2020)Mao, Jagtap, and Karniadakis}]{R26}
\bibinfo{author}{Z.~Mao}, \bibinfo{author}{A.~D. Jagtap},
  \bibinfo{author}{G.~E. Karniadakis},
\newblock \bibinfo{title}{Physics-informed neural networks for high-speed
  flows},
\newblock \bibinfo{journal}{Computer Methods in Applied Mechanics and
  Engineering} \bibinfo{volume}{360} (\bibinfo{year}{2020})
  \bibinfo{pages}{112789}.
\bibitem[{McClenny and Braga-Neto(2023)}]{R27}
\bibinfo{author}{L.~D. McClenny}, \bibinfo{author}{U.~M. Braga-Neto},
\newblock \bibinfo{title}{Self-adaptive physics-informed neural networks},
\newblock \bibinfo{journal}{Journal of Computational Physics}
  \bibinfo{volume}{474} (\bibinfo{year}{2023}) \bibinfo{pages}{111722}.
\bibitem[{Henkes et~al.(2022)Henkes, Wessels, and Mahnken}]{R29}
\bibinfo{author}{A.~Henkes}, \bibinfo{author}{H.~Wessels},
  \bibinfo{author}{R.~Mahnken},
\newblock \bibinfo{title}{Physics informed neural networks for continuum
  micromechanics},
\newblock \bibinfo{journal}{Computer Methods in Applied Mechanics and
  Engineering} \bibinfo{volume}{393} (\bibinfo{year}{2022})
  \bibinfo{pages}{114790}.
\bibitem[{Yu et~al.(2018)}]{R18}
\bibinfo{author}{B.~Yu}, et~al.,
\newblock \bibinfo{title}{The deep ritz method: a deep learning-based numerical
  algorithm for solving variational problems},
\newblock \bibinfo{journal}{Communications in Mathematics and Statistics}
  \bibinfo{volume}{6} (\bibinfo{year}{2018}) \bibinfo{pages}{1--12}.
\bibitem[{Sirignano and Spiliopoulos(2018)}]{R19}
\bibinfo{author}{J.~Sirignano}, \bibinfo{author}{K.~Spiliopoulos},
\newblock \bibinfo{title}{Dgm: A deep learning algorithm for solving partial
  differential equations},
\newblock \bibinfo{journal}{Journal of computational physics}
  \bibinfo{volume}{375} (\bibinfo{year}{2018}) \bibinfo{pages}{1339--1364}.
\bibitem[{Samaniego et~al.(2020)Samaniego, Anitescu, Goswami, Nguyen-Thanh,
  Guo, Hamdia, Zhuang, and Rabczuk}]{R20}
\bibinfo{author}{E.~Samaniego}, \bibinfo{author}{C.~Anitescu},
  \bibinfo{author}{S.~Goswami}, \bibinfo{author}{V.~M. Nguyen-Thanh},
  \bibinfo{author}{H.~Guo}, \bibinfo{author}{K.~Hamdia},
  \bibinfo{author}{X.~Zhuang}, \bibinfo{author}{T.~Rabczuk},
\newblock \bibinfo{title}{An energy approach to the solution of partial
  differential equations in computational mechanics via machine learning:
  Concepts, implementation and applications},
\newblock \bibinfo{journal}{Computer Methods in Applied Mechanics and
  Engineering} \bibinfo{volume}{362} (\bibinfo{year}{2020})
  \bibinfo{pages}{112790}.
\bibitem[{Peng et~al.(2023)Peng, Hu, and Xu}]{R28}
\bibinfo{author}{Y.~Peng}, \bibinfo{author}{D.~Hu}, \bibinfo{author}{Z.-Q.~J.
  Xu},
\newblock \bibinfo{title}{A non-gradient method for solving elliptic partial
  differential equations with deep neural networks},
\newblock \bibinfo{journal}{Journal of Computational Physics}
  \bibinfo{volume}{472} (\bibinfo{year}{2023}) \bibinfo{pages}{111690}.
\bibitem[{Lyu et~al.(2022)Lyu, Zhang, Chen, and Chen}]{R36}
\bibinfo{author}{L.~Lyu}, \bibinfo{author}{Z.~Zhang},
  \bibinfo{author}{M.~Chen}, \bibinfo{author}{J.~Chen},
\newblock \bibinfo{title}{Mim: A deep mixed residual method for solving
  high-order partial differential equations},
\newblock \bibinfo{journal}{Journal of Computational Physics}
  \bibinfo{volume}{452} (\bibinfo{year}{2022}) \bibinfo{pages}{110930}.
\bibitem[{Uriarte et~al.(2023)Uriarte, Pardo, Muga, and Mu{\~n}oz-Matute}]{R48}
\bibinfo{author}{C.~Uriarte}, \bibinfo{author}{D.~Pardo},
  \bibinfo{author}{I.~Muga}, \bibinfo{author}{J.~Mu{\~n}oz-Matute},
\newblock \bibinfo{title}{A deep double ritz method (d2rm) for solving partial
  differential equations using neural networks},
\newblock \bibinfo{journal}{Computer Methods in Applied Mechanics and
  Engineering} \bibinfo{volume}{405} (\bibinfo{year}{2023})
  \bibinfo{pages}{115892}.
\bibitem[{Yang et~al.(2022)Yang, Zuo, Tian, and Lei}]{R51}
\bibinfo{author}{J.~Yang}, \bibinfo{author}{J.~Zuo}, \bibinfo{author}{Y.~Tian},
  \bibinfo{author}{M.~Lei},
\newblock \bibinfo{title}{Deep ritz method for solving high-dimensional
  fractional differential equations},
\newblock \bibinfo{journal}{Engineered Science} \bibinfo{volume}{21}
  (\bibinfo{year}{2022}) \bibinfo{pages}{789}.
\bibitem[{Lai et~al.(2022)Lai, Chang, Lin, Hu, and Lin}]{R52}
\bibinfo{author}{M.-C. Lai}, \bibinfo{author}{C.-C. Chang},
  \bibinfo{author}{W.-S. Lin}, \bibinfo{author}{W.-F. Hu},
  \bibinfo{author}{T.-S. Lin},
\newblock \bibinfo{title}{A shallow ritz method for elliptic problems with
  singular sources},
\newblock \bibinfo{journal}{Journal of Computational Physics}
  \bibinfo{volume}{469} (\bibinfo{year}{2022}) \bibinfo{pages}{111547}.
\bibitem[{Wang and Zhang(2020)}]{R53}
\bibinfo{author}{Z.~Wang}, \bibinfo{author}{Z.~Zhang},
\newblock \bibinfo{title}{A mesh-free method for interface problems using the
  deep learning approach},
\newblock \bibinfo{journal}{Journal of Computational Physics}
  \bibinfo{volume}{400} (\bibinfo{year}{2020}) \bibinfo{pages}{108963}.
\bibitem[{Minakowski and Richter(2023)}]{R54}
\bibinfo{author}{P.~Minakowski}, \bibinfo{author}{T.~Richter},
\newblock \bibinfo{title}{A priori and a posteriori error estimates for the
  deep ritz method applied to the laplace and stokes problem},
\newblock \bibinfo{journal}{Journal of Computational and Applied Mathematics}
  \bibinfo{volume}{421} (\bibinfo{year}{2023}) \bibinfo{pages}{114845}.
\bibitem[{Leung et~al.(2022)Leung, Lin, and Zhang}]{R32}
\bibinfo{author}{W.~T. Leung}, \bibinfo{author}{G.~Lin},
  \bibinfo{author}{Z.~Zhang},
\newblock \bibinfo{title}{Nh-pinn: Neural homogenization-based physics-informed
  neural network for multiscale problems},
\newblock \bibinfo{journal}{Journal of Computational Physics}
  \bibinfo{volume}{470} (\bibinfo{year}{2022}) \bibinfo{pages}{111539}.
\bibitem[{Liu et~al.(2020)Liu, Cai, and Xu}]{R30}
\bibinfo{author}{Z.~Liu}, \bibinfo{author}{W.~Cai}, \bibinfo{author}{Z.-Q.~J.
  Xu},
\newblock \bibinfo{title}{Multi-scale deep neural network (mscalednn) for
  solving poisson-boltzmann equation in complex domains},
\newblock \bibinfo{journal}{arXiv preprint arXiv:2007.11207}
  (\bibinfo{year}{2020}).
\bibitem[{Li et~al.(2020)Li, Xu, and Zhang}]{R31}
\bibinfo{author}{X.-A. Li}, \bibinfo{author}{Z.-Q.~J. Xu},
  \bibinfo{author}{L.~Zhang},
\newblock \bibinfo{title}{A multi-scale dnn algorithm for nonlinear elliptic
  equations with multiple scales},
\newblock \bibinfo{journal}{arXiv preprint arXiv:2009.14597}
  (\bibinfo{year}{2020}).
\bibitem[{Li et~al.(2023)Li, Xu, and Zhang}]{R33}
\bibinfo{author}{X.-A. Li}, \bibinfo{author}{Z.-Q.~J. Xu},
  \bibinfo{author}{L.~Zhang},
\newblock \bibinfo{title}{Subspace decomposition based dnn algorithm for
  elliptic type multi-scale pdes},
\newblock \bibinfo{journal}{Journal of Computational Physics}
  (\bibinfo{year}{2023}) \bibinfo{pages}{112242}.
\bibitem[{Jiang et~al.(2023)Jiang, Wu, Chen, Chatzigeorgiou, and
  Meraghni}]{R37}
\bibinfo{author}{J.~Jiang}, \bibinfo{author}{J.~Wu}, \bibinfo{author}{Q.~Chen},
  \bibinfo{author}{G.~Chatzigeorgiou}, \bibinfo{author}{F.~Meraghni},
\newblock \bibinfo{title}{Physically informed deep homogenization neural
  network for unidirectional multiphase/multi-inclusion thermoconductive
  composites},
\newblock \bibinfo{journal}{Computer Methods in Applied Mechanics and
  Engineering} \bibinfo{volume}{409} (\bibinfo{year}{2023})
  \bibinfo{pages}{115972}.
\bibitem[{Xu and Zhou(2021)}]{R49}
\bibinfo{author}{Z.~J. Xu}, \bibinfo{author}{H.~Zhou},
\newblock \bibinfo{title}{Deep frequency principle towards understanding why
  deeper learning is faster},
\newblock in: \bibinfo{booktitle}{Proceedings of the AAAI Conference on
  Artificial Intelligence}, volume~\bibinfo{volume}{35}, \bibinfo{year}{2021},
  pp. \bibinfo{pages}{10541--10550}.
\bibitem[{Xu et~al.(2019)Xu, Zhang, and Xiao}]{R50}
\bibinfo{author}{Z.-Q.~J. Xu}, \bibinfo{author}{Y.~Zhang},
  \bibinfo{author}{Y.~Xiao},
\newblock \bibinfo{title}{Training behavior of deep neural network in frequency
  domain},
\newblock in: \bibinfo{booktitle}{Neural Information Processing: 26th
  International Conference, ICONIP 2019, Sydney, NSW, Australia, December
  12--15, 2019, Proceedings, Part I 26}, \bibinfo{organization}{Springer},
  \bibinfo{year}{2019}, pp. \bibinfo{pages}{264--274}.
\bibitem[{Dong et~al.(2018)Dong, Cui, Nie, Yang, and Yang}]{R42}
\bibinfo{author}{H.~Dong}, \bibinfo{author}{J.~Cui}, \bibinfo{author}{Y.~Nie},
  \bibinfo{author}{Z.~Yang}, \bibinfo{author}{Z.~Yang},
\newblock \bibinfo{title}{Multiscale computational method for heat conduction
  problems of composite structures with diverse periodic configurations in
  different subdomains},
\newblock \bibinfo{journal}{Computers \& Mathematics with Applications}
  \bibinfo{volume}{76} (\bibinfo{year}{2018}) \bibinfo{pages}{2549--2565}.
\bibitem[{Wang et~al.(2015)Wang, Cao, and Wong}]{R43}
\bibinfo{author}{X.~Wang}, \bibinfo{author}{L.~Cao}, \bibinfo{author}{Y.~Wong},
\newblock \bibinfo{title}{Multiscale computation and convergence for coupled
  thermoelastic system in composite materials},
\newblock \bibinfo{journal}{Multiscale Modeling \& Simulation}
  \bibinfo{volume}{13} (\bibinfo{year}{2015}) \bibinfo{pages}{661--690}.
\bibitem[{Cao(2006)}]{R44}
\bibinfo{author}{L.-Q. Cao},
\newblock \bibinfo{title}{Multiscale asymptotic expansion and finite element
  methods for the mixed boundary value problems of second order elliptic
  equation in perforated domains},
\newblock \bibinfo{journal}{Numerische Mathematik} \bibinfo{volume}{103}
  (\bibinfo{year}{2006}) \bibinfo{pages}{11--45}.
\bibitem[{He et~al.(2016)He, Zhang, Ren, and Sun}]{R35}
\bibinfo{author}{K.~He}, \bibinfo{author}{X.~Zhang}, \bibinfo{author}{S.~Ren},
  \bibinfo{author}{J.~Sun},
\newblock \bibinfo{title}{Deep residual learning for image recognition},
\newblock in: \bibinfo{booktitle}{Proceedings of the IEEE conference on
  computer vision and pattern recognition}, \bibinfo{year}{2016}, pp.
  \bibinfo{pages}{770--778}.
\bibitem[{Jiao et~al.(2021)Jiao, Lai, Li, Lu, Wang, Wang, and Yang}]{R34}
\bibinfo{author}{Y.~Jiao}, \bibinfo{author}{Y.~Lai}, \bibinfo{author}{D.~Li},
  \bibinfo{author}{X.~Lu}, \bibinfo{author}{F.~Wang},
  \bibinfo{author}{Y.~Wang}, \bibinfo{author}{J.~Z. Yang},
\newblock \bibinfo{title}{A rate of convergence of physics informed neural
  networks for the linear second order elliptic pdes},
\newblock \bibinfo{journal}{arXiv preprint arXiv:2109.01780}
  (\bibinfo{year}{2021}).
\bibitem[{Lu et~al.(2021)Lu, Pestourie, Yao, Wang, Verdugo, and Johnson}]{R55}
\bibinfo{author}{L.~Lu}, \bibinfo{author}{R.~Pestourie},
  \bibinfo{author}{W.~Yao}, \bibinfo{author}{Z.~Wang},
  \bibinfo{author}{F.~Verdugo}, \bibinfo{author}{S.~G. Johnson},
\newblock \bibinfo{title}{Physics-informed neural networks with hard
  constraints for inverse design},
\newblock \bibinfo{journal}{SIAM Journal on Scientific Computing}
  \bibinfo{volume}{43} (\bibinfo{year}{2021}) \bibinfo{pages}{B1105--B1132}.
\bibitem[{Hornik(1991)}]{R56}
\bibinfo{author}{K.~Hornik},
\newblock \bibinfo{title}{Approximation capabilities of multilayer feedforward
  networks},
\newblock \bibinfo{journal}{Neural networks} \bibinfo{volume}{4}
  (\bibinfo{year}{1991}) \bibinfo{pages}{251--257}.

\end{thebibliography}







\end{document}